\documentclass[final,10pt,3p]{elsarticle}
\journal{NeuroImage}
\usepackage[utf8]{inputenc}
\usepackage[T1]{fontenc}
\usepackage{newunicodechar}
\newunicodechar{−}{\ensuremath{-}}
\newunicodechar{≤}{\ensuremath{\leq}}
\newunicodechar{≥}{\ensuremath{\geq}}
\usepackage[english]{babel}
\usepackage{amsmath,amssymb,amsthm,mathtools}
\usepackage{graphicx}
\usepackage{float}
\usepackage{subcaption}
\usepackage[labelfont=bf]{caption}
\usepackage{booktabs}
\usepackage{tabularx}
\usepackage{multirow}
\usepackage{colortbl}
\usepackage{siunitx}
\usepackage{hyperref}
\hypersetup{colorlinks=true,linkcolor=blue,citecolor=blue,urlcolor=blue}
\usepackage[noabbrev,capitalize,nameinlink]{cleveref}
\usepackage{enumitem}
\usepackage{pdflscape}
\usepackage{changepage}
\DeclareOldFontCommand{\bf}{\normalfont\bfseries}{\mathbf}
\makeatletter\renewcommand{\fnum@figure}{Fig. \thefigure}\makeatother
\Crefname{figure}{Fig.}{Figs.}
\usepackage{xcolor}

\journal{NeuroImage}
\date{}
\title{}
\hypersetup{
 pdfauthor={},
 pdftitle={},
 pdfkeywords={},
 pdfsubject={},
 pdfcreator={},
 pdflang={English}}
\usepackage{natbib}
\begin{document}

\begin{frontmatter}

\title{Homology-based Morphometry of Brain Atrophy: Methods and Applications}

\author[math]{D. Quiccione}
\author[math]{M. Pirashvili}
\author[math]{N. Broomhead}
\author[psy]{S. J. Fallon}

\affiliation[math]{organization={School of Engineering, Computing and Mathematics, University of Plymouth},
            city={Plymouth},
            country={United Kingdom}}

\affiliation[psy]{organization={School of Psychology, University of Plymouth},
            city={Plymouth},
            country={United Kingdom}}

\begin{abstract}
Understanding the structure of the brain, and how it changes with time and disease, is a core goal of structural neuroimaging. Contemporary approaches to structural brain analysis are dominated by voxel-wise, mass-univariate methods such as voxel-based morphometry (VBM). However, these techniques require images to be normalized to a standard template, which can obscure subject-specific geometric features. Normalization to a common stereotactic space can also be problematic when comparing groups with substantial brain pathology, lesions, or other anatomical abnormalities.

Here, we introduce two complementary pipelines based on persistent homology (PH), a tool from topological data analysis, to quantify multiscale geometric features of structural T1-weighted MRI scans. Pipeline 1 quantifies regional thinning by applying the Euclidean distance transform to tissue masks in a slice-wise manner. Pipeline 2 uses \(\alpha\)-filtrations to measure structural similarity between pairs of scans, capturing sulcal widening and ventricular enlargement.

Synthetic experiments with controlled induced lesions showed that Pipeline 1 is best suited to between-subject analyses, whereas Pipeline 2 is better suited to within-subject designs. Applied to real-world data from the Alzheimer's Disease Neuroimaging Initiative (ADNI), Pipeline 1 separated Alzheimer's disease (AD) from cognitively normal (CN) participants using single-modality T1-weighted MRI without nonlinear registration (ROC-AUC = 0.895), with peak effects localized to medial temporal regions. Pipeline 2 captured disease-related longitudinal change, with follow-up scans remaining closest to their own baselines and AD subjects showing greater short-interval change than CN subjects. Together, these pipelines provide interpretable topological biomarkers for cross-sectional group comparisons and longitudinal tracking.


\end{abstract}



\end{frontmatter}
\section{Introduction}
\label{sec:orge277aab}

A key goal of human neuroscience, both basic and clinical, is to characterize the structure and shape of the living brain \citep{Sebenius_2024}. One of the most widely used tools for this purpose is magnetic resonance imaging (MRI), particularly T1-weighted (T1w) images, which form the mainstay of structural brain imaging because of their safety and high spatial resolution. The utility of T1w MRI has been greatly expanded by automated analysis packages such as voxel-based morphometry (VBM) \citep{Ashburner2000VBM}. VBM has yielded fundamental insights into the relationship between brain and behaviour \citep{Genon_2022} and neurodevelopmental and neurodegenerative disorders \citep{Liloia_2024,Ferreira_2011,Rowe_2010}.

\begin{figure*}[h]
\centering
\includegraphics[width=0.5\textwidth]{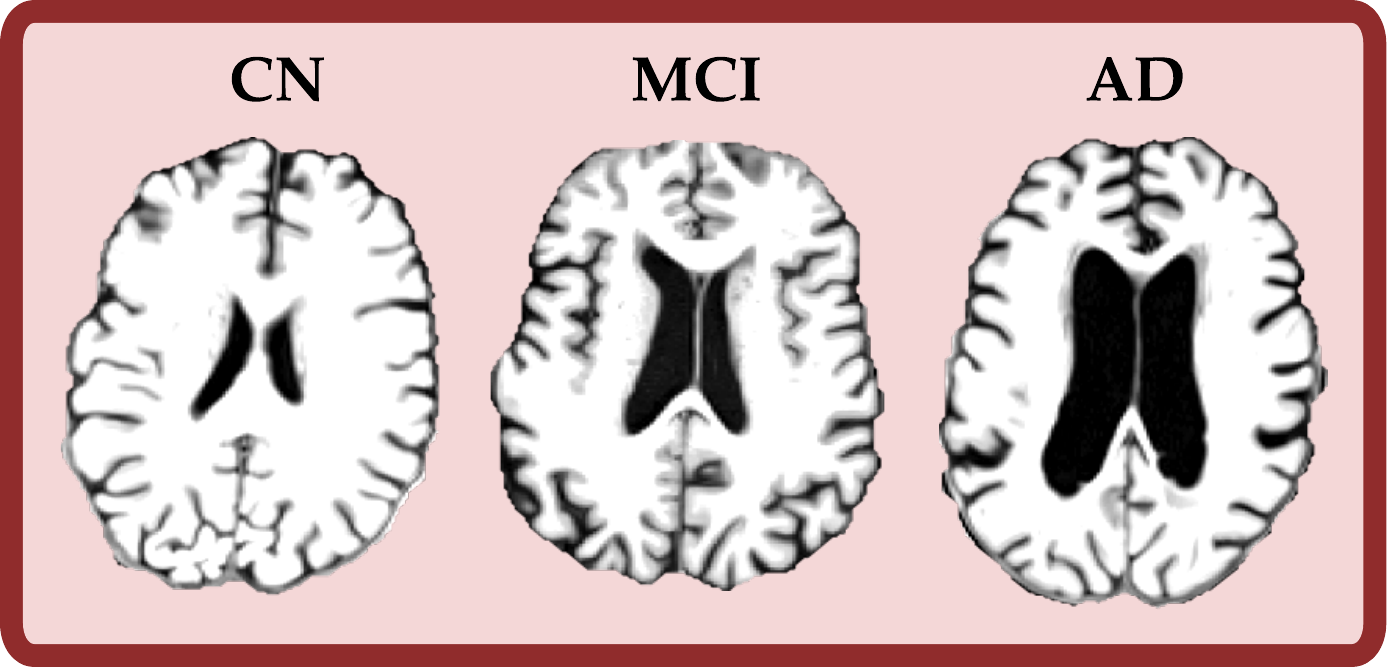}
\caption{\label{fig:cn-mci-ad-panel}Brain atrophy progression visualized in T1w MRI. On the left, preserved tissue in a cognitively normal (CN) participant. In the middle, intermediate atrophy with early ventricular expansion in a mild cognitive impairment (MCI) participant. On the right, advanced atrophy with marked ventricular dilation, sulcal widening, and cortical thinning in an AD participant. Data from ADNI.}
\end{figure*}

However, current morphometric frameworks still characterize brain change incompletely. ROI-based volumetry requires complex segmentations and atlases and reduces rich spatial organization to regional scalar summaries \citep{Yaakub_2020}. VBM offers better localization, but it relies on nonlinear registration to standard templates, introducing biases and reproducibility concerns \citep{Diaz_de_Grenu_2013,Zhou_2022,Giuliani_2005}. Its mass-univariate design also requires multiple-comparison correction, which remains problematic \citep{Henley_2009,Eklund_2016}. These approaches can therefore struggle to produce robust subject-level descriptors \citep{Scarpazza_2013}. Deep learning methods often achieve strong classification performance, but they commonly sacrifice interpretability and can generalize poorly across datasets \citep{Wen2020,Ebrahimighahnavieh2020}. Together, these limitations reduce sensitivity to the structural changes needed for diagnosis and disease tracking.

An emerging framework for studying the shape of data is topological data analysis (TDA). TDA has been used extensively in medical imaging and has successfully revealed patterns that are obscured by conventional methods \citep{Singh2023,Salch_2021,De_Benedictis_2024,Olave_2025}. One of its most powerful tools is persistent homology (PH), which analyzes how topological features of data---such as connected components (\(H_0\)), loops (\(H_1\)), and voids (\(H_2\))---appear and disappear across parameter scales \citep{CohenSteiner2007Stability,EdelsbrunnerHarer2010}.

Prior work has pioneered the application of PH to brain morphometry, particularly in brain tumors \citep{Francois2024,Crawford2020} and in structural connectivity and brain surface data \citep{Chung_2009,Chung_2015,Chung_2024,Phillips_2024}.

In neurodegeneration, Pachauri et al. \citep{Pachauri2011} compared the topology of cortical thickness maps in Alzheimer's disease with controls and reported good group separation (ROC-AUC 0.80), but their approach required computationally intensive surface reconstructions and complex preprocessing. PH has also been applied to hippocampal shape reconstructions, where coordinate-free descriptors reduced preprocessing demands for statistical inference \citep{Hofer2017}, although that study included only 42 participants. At the whole-brain level, Saadat-Yazdi et al. \citep{SaadatYazdi2021} computed topological summaries from gray- and white-matter masks but achieved only modest balanced accuracy (\(0.76 \pm 0.02\)). PH has also been used as input to convolutional neural networks for Alzheimer's disease diagnosis \citep{Bruningk2021}, but intensity-based filtrations of that kind lack direct pathological specificity.

Beyond these anatomical applications, most neuroimaging uses of PH have focused on functional, structural, and metabolic brain networks derived from atlas-dependent parcellations \citep{Glozman_2016,Salch_2021,Singh_2025}. What remains missing is a whole-brain PH framework that works directly in subject space, scales to clinically heterogeneous datasets, and avoids nonlinear normalization to a standard template.

Here, we introduce two complementary pipelines that address this gap. Both operate on standard T1w-derived tissue masks in subject space and require only routine tissue segmentation. Conceptually, they target two measurable consequences of brain atrophy: tissue thinning and compensatory CSF expansion.

\emph{Pipeline 1} captures tissue thinning through slice-wise Euclidean distance transform (EDT) superlevel filtrations of the parenchymal mask, computing \(H_1\) persistence to generate features for cross-sectional analysis. \emph{Pipeline 2} captures CSF expansion through \(\alpha\)-complex filtrations of the CSF complement mask, computing \(H_2\) persistence for within-subject longitudinal tracking. Thinning of brain tissue reduces EDT depths and shortens \(H_1\) cycle lifetimes, whereas expansion of CSF spaces increases \(H_2\) persistence. These descriptors therefore map directly onto biologically meaningful consequences of atrophy while preserving spatial organization beyond volumetric summaries.

We evaluated both pipelines in synthetic and clinical settings. In synthetic validation, both responded in the expected directions to controlled tissue erosion. In ADNI data, Pipeline 1 supported cross-sectional group separation (ROC-AUC = \(0.895\); AD vs CN) and anatomical localization, whereas Pipeline 2 supported within-subject longitudinal tracking. Topological features also correlated with hippocampal volume and cognitive measures, outperformed CSF volumetrics, and were computed within a fully reproducible workflow from raw DICOMs through final statistical inference.

Together, these pipelines provide a PH-based morphometric framework that detects and localizes brain atrophy, producing coordinate-free biomarkers in subject space from biologically motivated filtrations and minimal preprocessing.
\section{Methods}
\label{sec:org0741f5c}
Brain atrophy produces two direct morphological changes: tissue thinning and CSF expansion. Our framework is organized around these consequences. Pipeline 1 measures tissue thinning from the parenchymal mask, whereas Pipeline 2 measures CSF cavity expansion from the CSF complement mask. One key approach in this work is to employ binary segmentation masks rather than intensity-masked images. This reduces preprocessing overhead by avoiding the instability of intensity values and their sensitivity to the choice of normalization algorithm.
\subsection{Overview of the Pipelines}
\label{sec:org5805524}

\begin{figure*}
\centering
\includegraphics[width=0.9\textwidth]{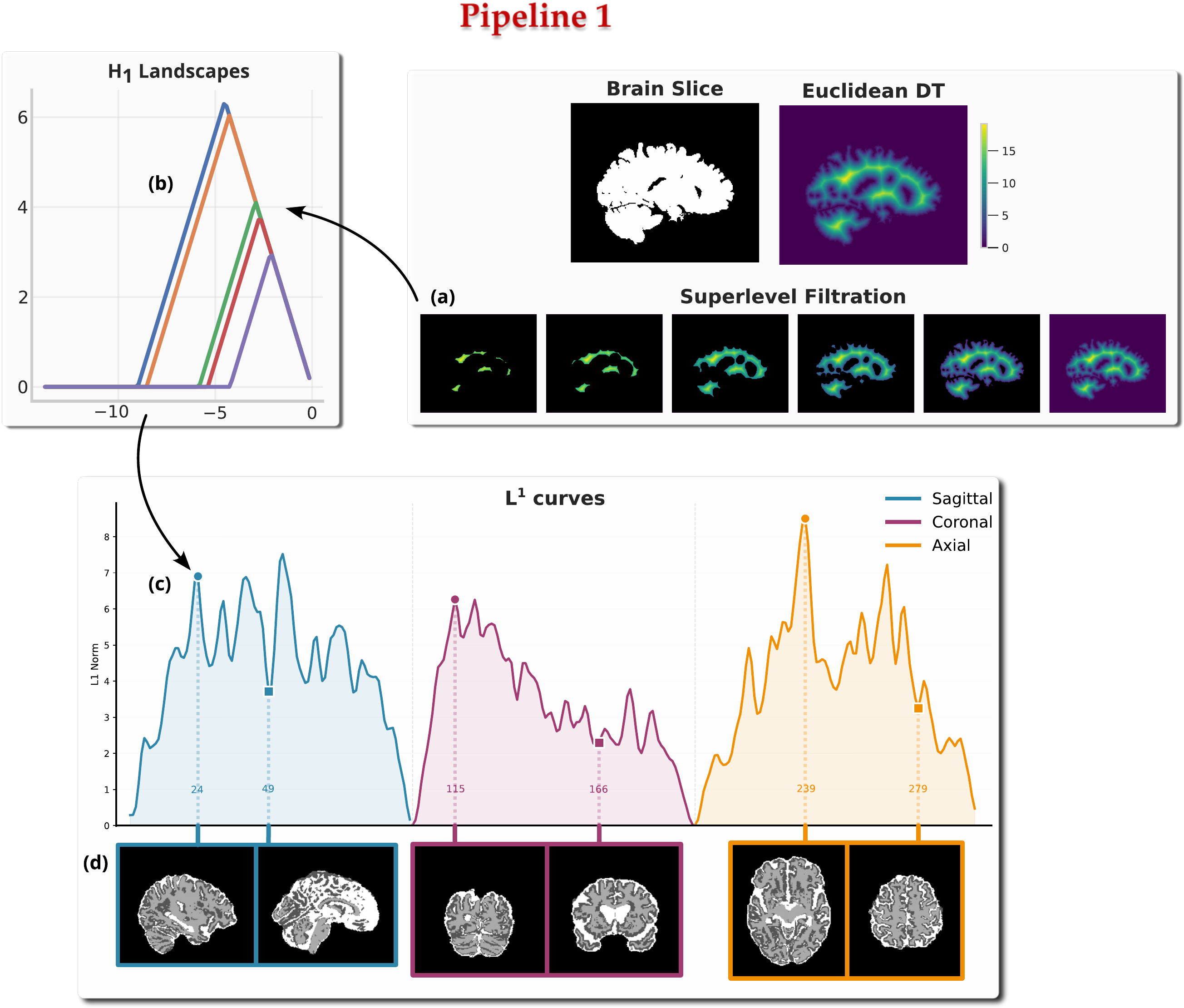}
\caption{\label{fig:pipeline1-info}Pipeline 1 schematic, (a) EDT superlevel filtration for one exemplar sagittal slice; (b) \(H_{1}\) landscapes associated to the filtration; (c) \(L^{1}\)-norm associated to the landscapes at the respective normalized slice position, and \(L^{1}\) curves in each anatomical view; (d) anatomical landmarks.}
\end{figure*}

\textbf{Pipeline 1} takes as input the \textbf{parenchymal mask}---the union of gray matter (GM) and white matter (WM).
Sets of parallel slices are taken in three directions, corresponding to the major anatomical views (sagittal, coronal and axial). For each slice, the EDT is computed---this assigns to each tissue voxel its distance from the nearest boundary, providing a slice-wise measure of local tissue depth.
The cubical \(H_1\) persistence of the EDT superlevel filtration is then calculated. Intuitively, this counts certain cycles within each slice whose depth is above a threshold \(p\), tracking when they appear and disappear as the value of \(p\) changes. These cycles represent holes, i.e., ring-like tissue configurations that encircle cavities or gaps (Fig. \ref{fig:pipeline1-info} (a)). Their persistence measures the erosion depth required to eliminate the hole. As neurodegeneration progresses, cortical thinning and subcortical volume loss reduce depth values, and this manifests as shortened lifetimes of cycles in our framework. 

Finally, the topological information is summarized into a single scalar value for each slice. This is done by converting the \(H_{1}\) persistent diagram for each slice into a persistence landscape and then taking its \(L^{1}\) norm (Bubenik \citep{Bubenik2015Landscapes}) (Fig. \ref{fig:pipeline1-info} (b)). Treating this scalar as a function of slice index gives one curve per anatomical direction. By repeating this for sagittal, coronal, and axial views, we obtain three \(L^{1}\)-curves for each subject (Fig. \ref{fig:pipeline1-info} (c)). These curves encode how topological features vary across the brain in the three directions and show characteristic peaks and valleys at major anatomical landmarks such as ventricles, corpus callosum, hippocampus and the cerebellum (Fig. \ref{fig:pipeline1-info} (d)). This workflow provides localized, registration-free measurements that remain reliable for cross-subject comparisons without requiring nonlinear normalization to a template.

\textbf{Pipeline 2} analyzes the \textbf{CSF complement mask}, i.e., the negative space of the CSF mask within a bounding box of the intracranial volume. To quantify the shape of CSF bounded cavities we use the \(H_{2}\) persistent diagram of the \(\alpha\)-complex filtration \citep{EdelsbrunnerHarer2010,Edelsbrunner1994AlphaShapes}. The \(H_{2}\) cycles are the natural choice, since they represent in this data structure exactly those voids that are enclosed by the CSF boundaries. Their persistence directly models the compensatory enlargement and redistribution of CSF spaces following tissue loss (Fig. \ref{fig:pipeline2-info}).
This pipeline is best suited for within-subject comparisons since persistence diagrams are a sensitive measure of the distribution of CSF bounded cavities and fingerprint the individual CSF geometry. Furthermore, the bottleneck distance offers a robust and stable measure of change due to atrophy.

\begin{figure*}
\centering
\includegraphics[width=\textwidth]{figures/pipeline2_infograph_V2.pdf}
\caption{\label{fig:pipeline2-info}Pipeline 2 schematic, (a) illustrates the alpha complex filtration on a 2D coronal slice; (b) shows the \(H_{2}\) persistent diagrams of the 3D alpha complex filtration for representative participants for each diagnosis group (AD and CN). The \(H_{2}\) diagrams of AD participants consistently show cycles with higher persistence.}
\end{figure*}

\begin{figure*}
\centering
\includegraphics[width=0.5\textwidth]{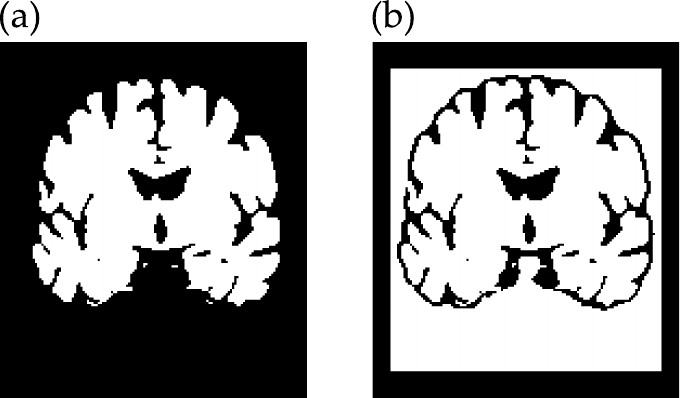}
\caption{\label{same-mask}2D coronal slices of the binary masks used as input to our pipelines. (a) Parenchymal mask \(P = GM \cup WM\) for Pipeline 1, in which we measure the local thickness of the tissue (white). (b) CSF complement mask for Pipeline 2, in which we measure the expansion of the CSF (black).}
\end{figure*}

Both pipelines start from the same tissue segmentation masks, and they differ mainly in whether we close the mask with a minimal bounding box (Fig. \ref{same-mask}).

Thus, our central hypotheses directly link topology to atrophy and its patterns: \((i)\) Parenchymal atrophy systematically reduces the area under the \(L^{1}\) curves; \((ii)\) CSF cavity expansion increases the persistence of \(H_2\) voids in \(\alpha\)-complex filtrations, yielding bottleneck distances between diagrams quantifying disease progression.

We dedicate a portion of the Experiments and Results section to test these hypotheses with controlled induced atrophy on CN subjects' parenchymal masks. 
\subsection{Computational implementations}
\label{sec:org33ca5e1}

All software components are open-source: \texttt{GUDHI} for alpha and cubical complex persistence computations \citep{Maria2014GUDHI,gudhi:urm}, \texttt{sMRIPrep} (v\(0.17.0\)) for tissue segmentation \citep{esteban_2024_14532379}, and \texttt{ANTs} for rigid alignment in Pipeline 1. 
\subsubsection{Pipeline 1.}
\label{sec:orga2cee25}
Each parenchymal mask \(P\) is resampled to \(1\,\text{mm}^{3}\) isotropic resolution and rigidly aligned (non-warping) to the MNI152 template to standardize orientation while preserving native geometry. For each anatomical view (sagittal, coronal, and axial), we extract 2D slices and compute the Euclidean distance transform (EDT), where each voxel stores its distance to the nearest boundary. Superlevel sets of EDT define a filtration that expands from deep tissue toward the boundary as the threshold decreases. Using 2D cubical homology \citep{gudhi:urm}, we compute \(H_1\) persistence diagrams per slice and convert them to persistence landscapes. The landscape \(L^1\) norm provides a scalar per slice; concatenating these across slices yields an anatomical \(L^{1}\) curve per view. Since EDT values are in mm, persistence coordinates (birth, death) and \(L^1\) norms have units of mm and mm\(^2\), respectively. Curves are interpolated to a common grid to enable comparison across subjects with varying brain sizes.
\subsubsection{Pipeline 2.}
\label{sec:orga6ad27c}

Each CSF complement mask is resampled to \(1\,\text{mm}^{3}\) isotropic resolution and cropped to the CSF bounding box. All voxel centers of non-CSF tissue form the input point cloud (dense sampling). The bounding box closure converts sulci from open grooves into enclosed cavities, so both ventricles and sulci contribute \(H_2\) features. The \(\alpha\)-complex filtration \citep{EdelsbrunnerHarer2010} grows simplices by squared circumradius, filling cavities as \(\alpha\) increases. Each \(H_2\) void appears at birth \(\alpha_b\) and disappears at death \(\alpha_d\) when filled. Coordinates are in mm\(^2\) (since \(\alpha = r^2\)); we apply a \(0.25\) mm\(^2\) persistence threshold to filter noise. Pairwise longitudinal comparisons use bottleneck distance \citep{EdelsbrunnerHarer2010} to quantify morphological change.

The stability of persistent homology \citep{CohenSteiner2007Stability} guarantees that small changes in segmentation result in small changes in the diagram. A bottleneck distance of \(\varepsilon\) mm\(^2\) represents the minimal perturbation required to match diagrams, where each cavity's birth and death values shift by at most \(\varepsilon\) to establish correspondence. This supports within-subject tracking: in our setting, increasing distances reflect cavity expansion over time.
\subsection{Data and Software}
\label{sec:org7aa884e}

This work uses data from the ADNI database \citep{Jack_2008} obtained from LONI Image and Data Archive (IDA). We analyze T1w structural MRI from subjects spanning the clinical spectrum from Cognitively Normal (CN) through (Early/Late) Mild Cognitive Impairment (EMCI, MCI and LMCI) to Alzheimer's disease (AD).

To minimize scanner-specific confounds, we restricted our analysis to T1-weighted MPRAGE (Magnetization Prepared Rapid Gradient Echo) sequences only. We obtained in this way two datasets (subsets of ADNI 1/2/GO), one for each pipeline.
For the cross-sectional dataset, we selected one visit per subject: the first available visit for CN and MCI subjects and the last available visit for AD subjects. This choice defines a validation setting intended to test whether the representation captures clear disease-related morphology. Age and TIV remained well-matched across groups after this selection (Table \ref{tab:demographics_3t}; one-way ANOVA: age \(F(2, 428) = 1.66\), \(p = 0.191\); TIV \(F(2, 428) = 2.56\), \(p = 0.078\)).
We applied Pipeline 1 to a cross-sectional \(3\)T dataset of \(n=765\) subjects (CN, EMCI, MCI, LMCI and AD). Binary classification (CN vs AD) uses \(n=334\) subjects. For the main results, we focused on a 3-group analysis (CN, MCI and AD) of \(n=431\) participants, to provide cleaner disease-stage contrasts. We report analyses of MCI subtypes (EMCI, LMCI) separately in the Supplement. Pipeline 2 analyzes a longitudinal (within-subject) \(1.5\)T dataset of \(n=82\) subjects, of which \(n=42\) CN and \(n=40\) AD. Covariate information is presented in the next section and the full table is presented in the Supplement.

We conducted our analysis starting from raw DICOM files to apply a consistent preprocessing pipeline to all images, as ADNI's preprocessed derivatives contained suboptimal or missing tissue segmentations.
The raw DICOM files were converted to NIfTI format and organized according to the Brain Imaging Data Structure (BIDS) specification using \texttt{dcm2bids} (v\(3.2.0\)). 
The structural images were preprocessed with \texttt{sMRIPrep}. Brain and tissue binary segmentation masks (GM, WM and CSF) were extracted from \texttt{dseg} derivatives. All segmentations were visually inspected for gross errors and misclassification of tissues (see Supplement).

For computational reproducibility, we controlled for known sources of stochasticity in neuroimaging pipelines. ANTs registration, while widely used, can introduce minor variability across runs due to random point sampling in metric computation and parallel processing. We mitigated these issues by setting a random seed for all registrations and enforcing single-threaded execution. These steps make feature extraction deterministic across runs.

The code for this project is fully available on GitLab \footnote{\url{https://gitlab.com/donatoqui/hbm}} for anyone to use and contribute to.
\section{Experiments and Results}
\label{sec:org2592dec}

\subsection*{Synthetic Erosion Alters Topology in Predicted Directions}
\label{sec:orga223737}
Before applying the pipelines to clinical data, we first tested whether they respond to controlled tissue loss in the predicted directions. To do this, we generated synthetic atrophy from cognitively normal ADNI participants, allowing us to isolate the effect of volume reduction without the heterogeneity of real disease. This experiment serves as a ground-truth validation of the biological logic introduced in the Methods: tissue thinning should reduce Pipeline 1's \(L^1\)-curve AUC, whereas the accompanying expansion of CSF-bounded cavities should increase Pipeline 2's \(H_2\) persistence and bottleneck distances.
\subsubsection*{Dataset Generation.}
\label{sec:org2b72a64}

\begin{figure}[htbp]
\centering
\includegraphics[width=0.65\textwidth]{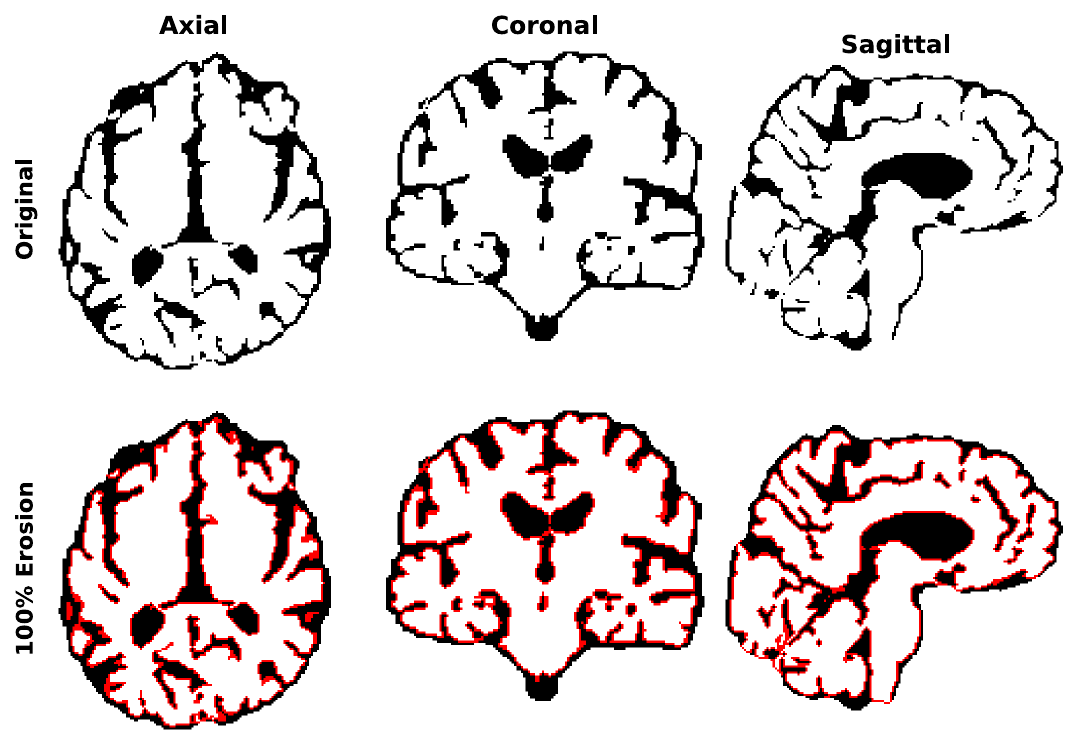}
\caption{\label{fig:synthetic_erosion}Example of synthetic erosion. The first row shows the original image; the second row shows the 100\% erosion, and in red the 1-voxel boundary layer around the parenchymal mask from which we stochastically removed \(25\%\), \(50\%\), \(75\%\), and \(100\%\) of voxels in the synthetic analysis in both pipelines.}
\end{figure}

We generated synthetic atrophy by applying stochastic partial erosion to cognitively normal (CN) parenchymal masks from ADNI. For each subject, we identified the 1-voxel boundary layer of the parenchymal mask (using a 6-connected structuring element) and stochastically removed \(0\%\), \(25\%\), \(50\%\), \(75\%\), and \(100\%\) of voxels from it (Fig. \ref{fig:synthetic_erosion}).

For Pipeline 1 (cross-sectional analysis), we generated \(200\) synthetic volumes from \(40\) 3\(T\) images at each erosion level. For Pipeline 2 (longitudinal analysis), we generated \(100\) synthetic volumes from \(20\) subjects acquired on \(1.5\)T scanners, matching the scanner field strength used in the subsequent clinical analyses.

\subsubsection*{L1-curve AUC Decreases with Tissue Loss.}
\label{sec:org0d0a13a}

\begin{figure*}
\centering
\includegraphics[width=0.99\textwidth]{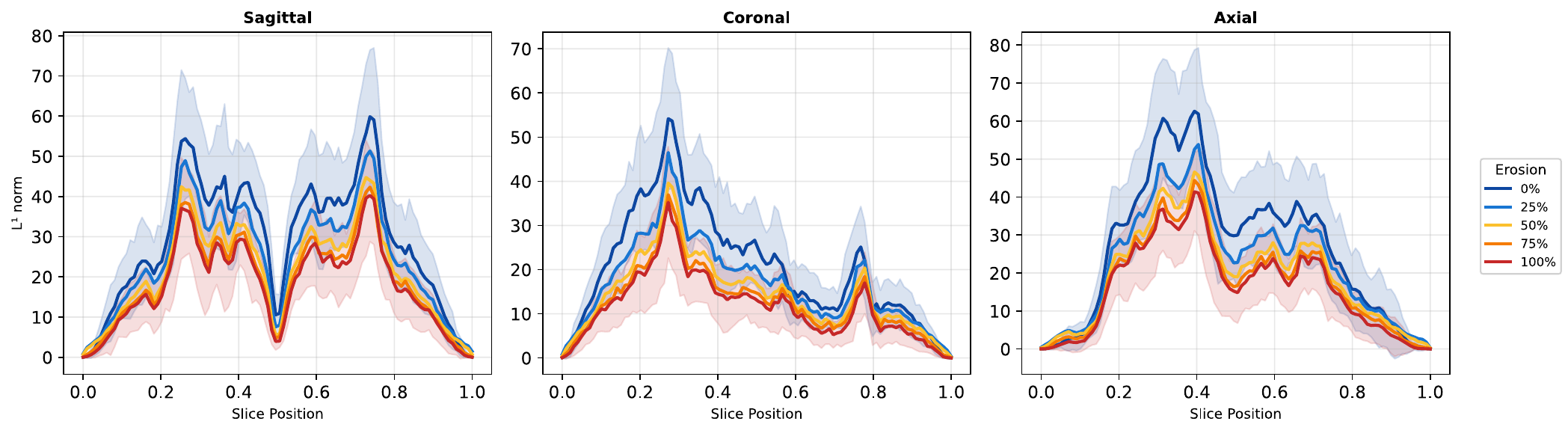}
\caption{\label{fig:synthetic_slices_validation}\(L^{1}\) curves for each synthetic erosion level, including standard deviation bands for the baseline (0\%) and 100\% erosion.}
\end{figure*}

The specifications above produced mean parenchymal volume reductions of \(4.6\%\), \(9.2\%\), \(13.8\%\), and \(18.4\%\) with respect to the \(25\%\), \(50\%\), \(75\%\), and \(100\%\) erosion levels.
Pipeline 1's \(L^{1}\) curves showed expected sensitivity to synthetic atrophy across all anatomical views (Fig. \ref{fig:synthetic_slices_validation}).

\(L^1\)-curve AUC decreased systematically with a \(38\%\), \(43\%\), and \(39\%\) reduction for sagittal, coronal, and axial views respectively for the \(100\%\) erosion level. This is consistent with AUC declining as tissue loss increases (Table \ref{tab:synthetic_correlation}). 

\begin{table*}[htbp]
\caption{\label{tab:synthetic_correlation}Correlation (Spearman rho) between volume loss and \(L^1\) curves AUC}
\centering
\small
\begin{tabular}{r c c}
\toprule
axis & Spearman rho & p value\\
Sagittal & -0.75 & 4.21e-37\\
Coronal & -0.76 & 8.20e-39\\
Axial & -0.76 & 6.31e-39\\
\bottomrule
\end{tabular}
\end{table*}
\subsubsection*{H2 Persistence Scales with Erosion.}
\label{sec:org23a8a2b}

\begin{figure*}[htbp]
\centering
\begin{tabular}{ll}
\includegraphics[width=.48\linewidth]{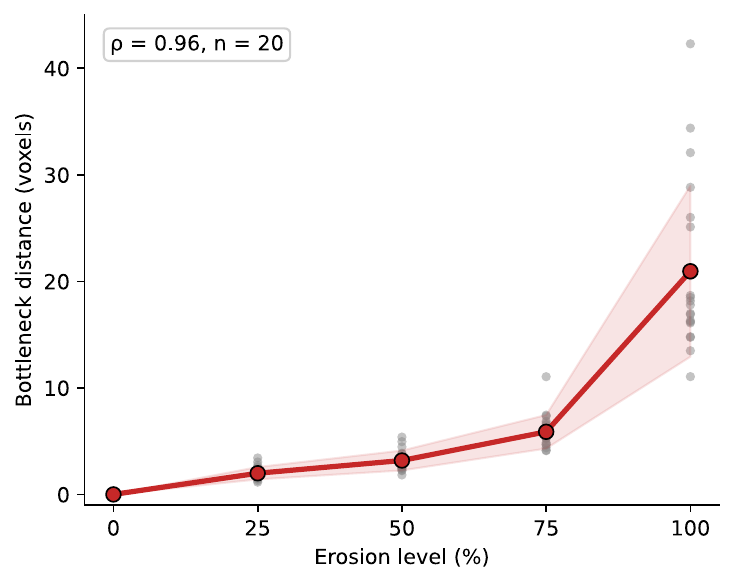}
 &
\includegraphics[width=.48\linewidth]{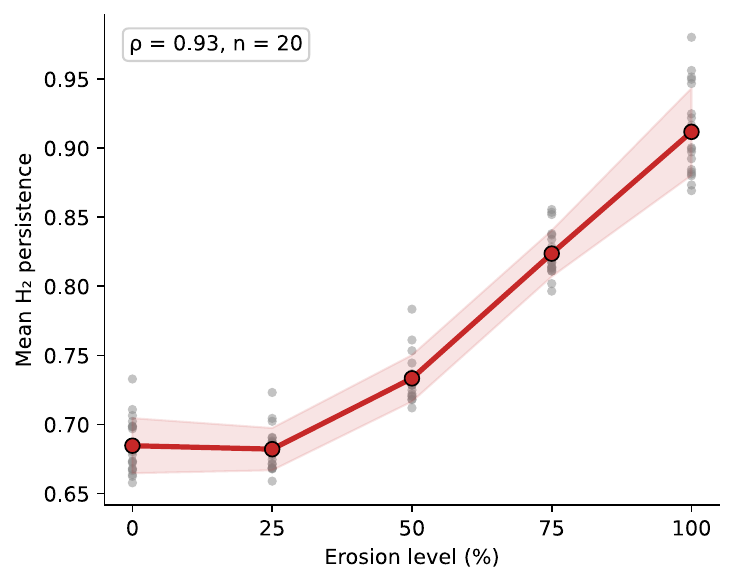}
\\
\end{tabular}
\caption{\label{fig:synthetic_alpha_validation} Synthetic validation for Pipeline 2. Left: \(H_2\) bottleneck distances from baseline (\(0\%\) erosion) increase with volume reduction. The points are individual subjects (\(n=20\)) at each erosion level, including the trivial baseline self-comparisons at \(0\%\) erosion. The red line is the mean trajectory, with SD shaded band. Spearman statistic excluding the baseline: \(\rho = 0.928\), \(p = 4.05 \times 10^{-35}\), \(n=80\). Right: mean finite \(H_2\) persistence increases across erosion levels.}
\end{figure*}

Pipeline 2 showed the expected persistence shift with synthetic cavity expansion in the \(20\)-subject dataset. From baseline to \(100\%\) erosion, mean \(H_2\) persistence increased from \(0.68\) to \(0.91\) mm\(^2\) (\(+33.8\%\)), and maximum \(H_2\) persistence increased from \(54.18\) to \(74.78\) mm\(^2\) (\(+38.0\%\)). Consistently, bottleneck distances from baseline increased with volume reduction (Fig. \ref{fig:synthetic_alpha_validation}); excluding the baseline, Spearman \(\rho = 0.928\), \(p = 4.05 \times 10^{-35}\), \(n=80\).

Both pipelines responded as expected to controlled tissue erosion, supporting sensitivity to tissue loss across the field strengths used in the subsequent clinical analyses (3T and 1.5T). In the next subsections we show that this sensitivity extends to real pathological data: Pipeline 1 detects and localizes cross-sectional group differences, and Pipeline 2 tracks longitudinal change.
\subsection*{Pipeline 1 Separates Disease Groups and Localizes Medial Temporal Effects}
\label{sec:org4abc71e}
\subsubsection*{Dataset demographics.}
\label{sec:org475dc18}

We apply Pipeline 1 to \(431\) subjects across three diagnostic groups: CN, MCI, and AD.

Binary classification analysis (CN vs AD) uses \(334\) subjects (201 CN + 133 AD), excluding MCI to focus on diagnostic extremes.

\begin{table*}[htbp]
\caption{\label{tab:demographics_3t}Demographics by diagnostic group (3T dataset). MMSE = Mini-Mental State Examination; TIV = total intracranial volume (mm\(^3\)).}
\centering
\small
\begin{tabular}{r c c c c}
\toprule
\textbf{Group} & \textbf{N (count)} & \textbf{Age (mean\(\pm\)SD)} & \textbf{MMSE (mean\(\pm\)SD)} & \textbf{TIV (mean\(\pm\)SD)}\\
\midrule
CN & 201 & 74.48\(\pm\) 6.70 & 28.94\(\pm\) 1.32 & \num{1424096}\(\pm\)\num{133255}\\
MCI & 97 & 75.22\(\pm\) 7.92 & 27.12\(\pm\) 2.31 & \num{1462362}\(\pm\)\num{129617}\\
AD & 133 & 75.98\(\pm\) 7.88 & 20.14\(\pm\) 4.86 & \num{1428998}\(\pm\)\num{157331}\\
\bottomrule
\end{tabular}
\end{table*}

Groups were well-matched on age and TIV (one-way ANOVA: age \(F(2, 428) = 1.66\), \(p = 0.191\); TIV \(F(2, 428) = 2.56\), \(p = 0.078\)).
\subsubsection*{L1 Curves Decrease with Disease Severity.}
\label{sec:org4baf598}

\(L^1\)-curve AUC decreases progressively from CN through MCI to AD, with similar anatomical patterns across subjects.

\begin{figure*}
\centering
\includegraphics[width=0.9\textwidth]{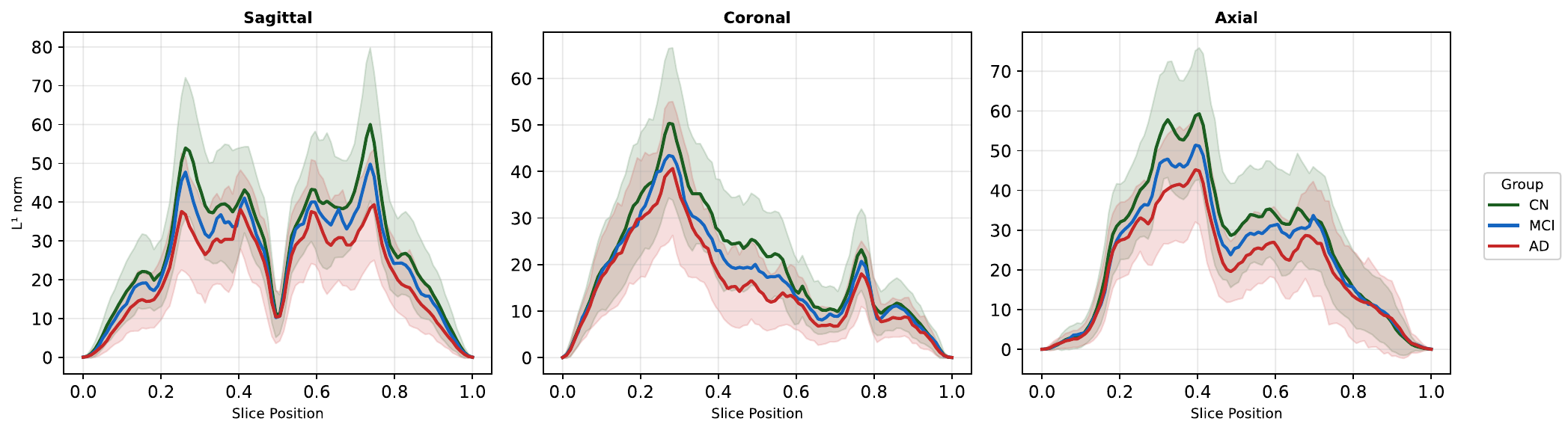}
\caption{\label{fig:l1_curves_all_groups}Mean \(L^1\) curves across diagnostic groups (CN, MCI, AD) for all three anatomical views. Shaded regions indicate SD for CN and AD groups. Progressive reduction in \(L^{1}\) norm from CN through MCI to AD tracks disease severity.}
\end{figure*}

The mean \(L^1\) curves (Fig. \ref{fig:l1_curves_all_groups}) separate progressively across diagnostic groups, with AD showing the lowest \(L^1\)-curve values across all three anatomical views and MCI occupying an intermediate position. Aggregated AUC features showed a large effect for CN vs AD (Cohen's \(d = 1.08\), \(p = 3.05 \times 10^{-20}\)) and moderate effects for adjacent stages: CN vs MCI (\(d = 0.47\), \(p = 3.26 \times 10^{-4}\)) and MCI vs AD (\(d = 0.57\), \(p = 6.27 \times 10^{-5}\)). This graded pattern is consistent with progressive tissue loss across disease severity rather than a binary separation between extremes. These group differences persisted after stratification by age and TIV (Supplement; Cohen's \(d\) range: \(0.51\)--\(1.51\)).

PCA of concatenated \(L^1\) curves further supports this interpretation, with disease severity emerging as the dominant mode of variation. The first two components explained \(40.31\%\) of total variance (PC1: \(33.30\%\), PC2: \(7.01\%\)). Further subgroup analyses (EMCI, LMCI) are provided in the Supplement.
\subsubsection*{Peak Effects Localize to Medial Temporal Regions.}
\label{sec:org816e95a}
We tested AD vs CN differences at each slice using Welch's \(t\)-tests with Bonferroni correction (\(n = 300\): \(100\) slices per axis \(\times\) \(3\) axes, \(\alpha = 0.05\)).

\begin{figure*}
\centering
\includegraphics[width=0.95\textwidth]{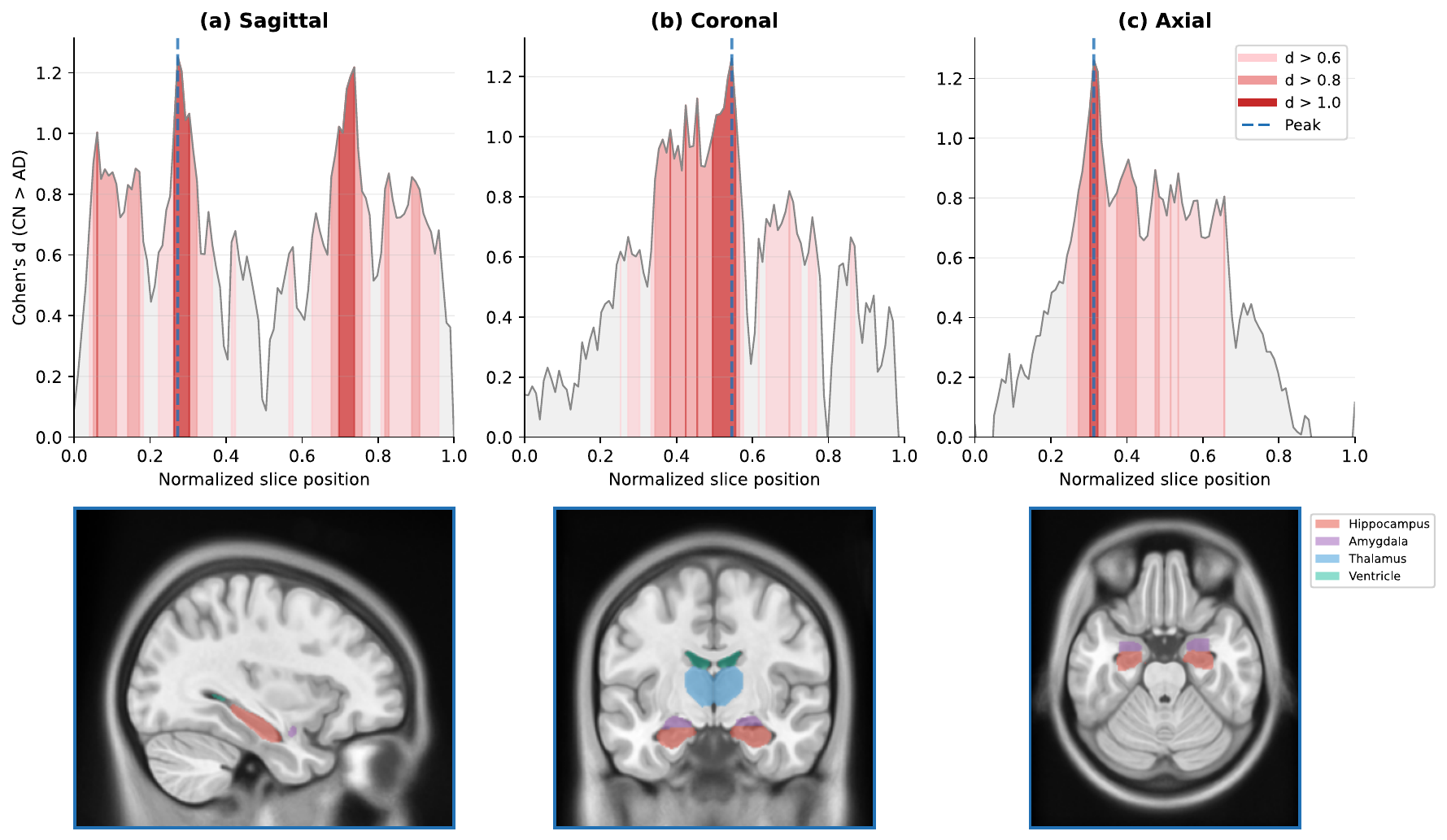}
\caption{\label{fig:slice_localization}Slice-wise effect size for (a) sagittal, (b) coronal, and (c) axial views. The dashed blue lines represent peaks at which the corresponding MNI slices (bottom) are shown.}
\end{figure*}

Peak effects clustered in medial temporal and connected limbic territories across all three views (Table \ref{tab:localization_pvalue}). The normalized position represents the relative position through each subject's brain along that axis, so position \(0.27\) corresponds to \(27\%\) of the way through the brain. The corresponding MNI slices intersect medial temporal and periventricular structures including the hippocampus, parahippocampal cortex, amygdala, thalamus, anterior cingulate cortex, and fusiform gyrus (Fig. \ref{fig:slice_localization}).

\begin{table}[htbp]
\caption{\label{tab:localization_pvalue}Peak localization of slice-wise group differences (Bonferroni-corrected \(p\)-values) and effect sizes (Cohen's \(d\), CN > AD)}
\centering
\small
\begin{tabular}{lcccc}
\toprule
Axis & Position & MNI (mm) & Cohen's \(d\) & \(p\)-value\\
\midrule
Sagittal & 0.27 & x = -33 & 1.25 & \(1.4 \times 10^{-24}\)\\
Coronal & 0.55 & y = -10 & 1.26 & \(4.1 \times 10^{-26}\)\\
Axial & 0.31 & z = -24 & 1.26 & \(1.8 \times 10^{-22}\)\\
\bottomrule
\end{tabular}
\end{table}
\subsubsection*{Combined Views Achieve Strong Classification.}
\label{sec:org5e8691b}

We evaluated classification performance using 5-fold stratified cross-validation with three models: PCA95+sLDA (\(95\%\) variance retention with shrinkage LDA), RBF-SVM (radial basis function support vector machine), and XGBoost. All preprocessing (scaling, PCA) was performed within scikit-learn Pipelines fit exclusively on training folds.

\begin{table*}[htbp]
\caption{\label{tab:classification}Classification performance metrics (5-fold cross-validation) for CN vs AD binary classification using landscape \(L^1\) curves}
\centering
\footnotesize
\begin{tabular}{llcccccc}
\toprule
Axis & Model & \multicolumn{3}{c}{Mean} & \multicolumn{3}{c}{Std Dev}\\
 &  & Bal. Acc. & F1 & ROC-AUC & Bal. Acc. & F1 & ROC-AUC\\
\midrule
All & \textbf{PCA95+sLDA} & \textbf{0.814} & \textbf{0.778} & \textbf{0.895} & \textbf{0.048} & \textbf{0.057} & \textbf{0.034}\\
 & RBF-SVM & 0.800 & 0.760 & 0.889 & 0.039 & 0.048 & 0.027\\
 & XGBoost & 0.775 & 0.728 & 0.874 & 0.046 & 0.055 & 0.029\\
\midrule
Axi & PCA95+sLDA & 0.732 & 0.665 & 0.824 & 0.035 & 0.050 & 0.034\\
 & RBF-SVM & 0.704 & 0.637 & 0.782 & 0.025 & 0.037 & 0.053\\
 & XGBoost & 0.779 & 0.729 & 0.835 & 0.030 & 0.040 & 0.035\\
\midrule
Cor & PCA95+sLDA & 0.728 & 0.674 & 0.845 & 0.040 & 0.048 & 0.034\\
 & RBF-SVM & 0.748 & 0.693 & 0.836 & 0.060 & 0.072 & 0.044\\
 & XGBoost & 0.755 & 0.704 & 0.853 & 0.051 & 0.063 & 0.048\\
\midrule
Sag & PCA95+sLDA & 0.764 & 0.719 & 0.865 & 0.042 & 0.043 & 0.035\\
 & RBF-SVM & 0.768 & 0.721 & 0.845 & 0.064 & 0.072 & 0.042\\
 & XGBoost & 0.762 & 0.714 & 0.855 & 0.037 & 0.042 & 0.047\\
\bottomrule
\end{tabular}
\end{table*}

Combining all three anatomical views achieved the best performance (PCA95+sLDA: balanced accuracy \(= 0.814\), ROC-AUC \(= 0.895\)). Individual anatomical views showed comparable performance, with sagittal sections providing the strongest single-view classification (balanced accuracy \(= 0.764\), ROC-AUC \(= 0.865\)). These results were obtained from single-modality T1w sMRI without additional imaging biomarkers.

Classification of adjacent disease stages proved more challenging: CN vs MCI achieved ROC-AUC \(= 0.699 \pm 0.028\) (XGBoost; balanced accuracy \(= 0.614 \pm 0.046\)) and MCI vs AD achieved ROC-AUC \(= 0.748 \pm 0.019\) (RBF-SVM; balanced accuracy \(= 0.667 \pm 0.043\)). Further MCI subtypes classification results are presented in the Supplement.
\subsection*{Pipeline 2 Yields Subject-Specific Longitudinal Fingerprints}
\label{sec:org874bcd3}
\subsubsection*{Dataset Characteristics.}
\label{sec:org9d16b38}
Pipeline 2 uses \(82\) subjects with paired visits separated by \(182 \pm 11\) days (Table \ref{tab:demographics_longitudinal}).

\begin{table}[htbp]
\caption{\label{tab:demographics_longitudinal}Longitudinal cohort characteristics by diagnostic group (baseline demographics)}
\centering
\small
\begin{tabular}{lcccc}
\toprule
Group & N & Age (mean\(\pm\)SD)\\
\midrule
AD & 40 & 76.0\(\pm\)8.0\\
CN & 42 & 77.6\(\pm\)4.4\\
\bottomrule
\end{tabular}
\end{table}
\subsubsection*{Follow-Up Scans Are Nearest to Own Baselines.}
\label{sec:org7652ddc}

\begin{figure*}[htbp]
\centering
\begin{tabular}{ll}
\includegraphics[width=.48\linewidth]{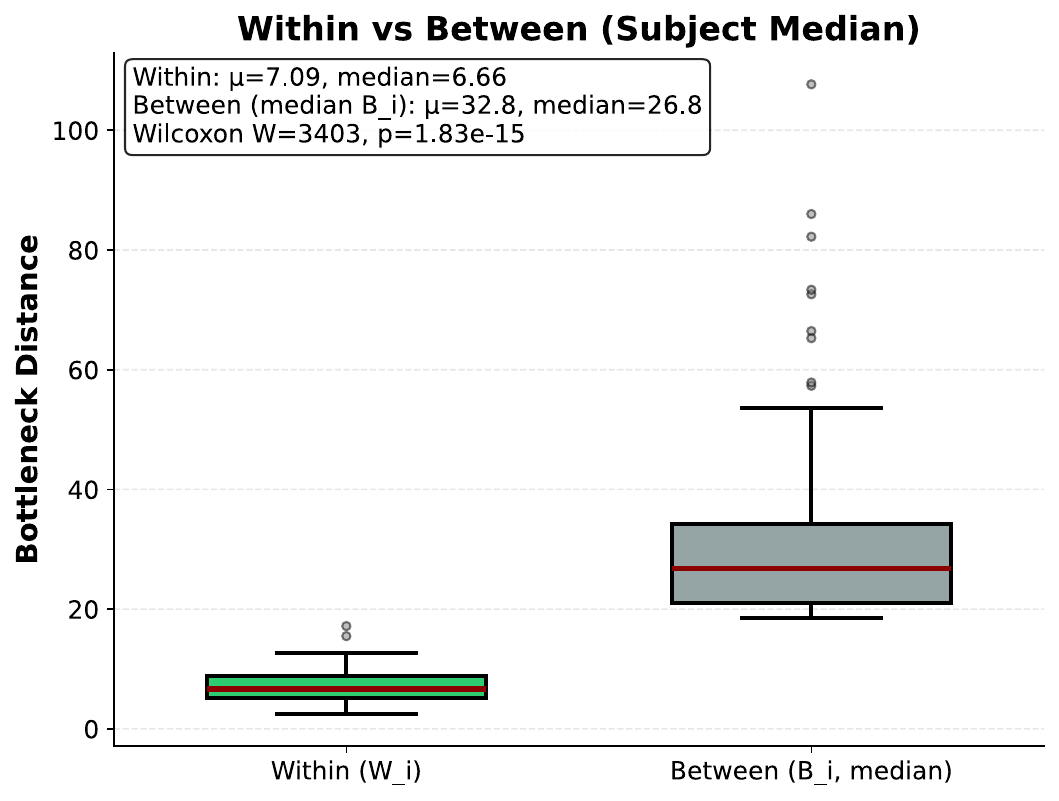}
 &
\includegraphics[width=.48\linewidth]{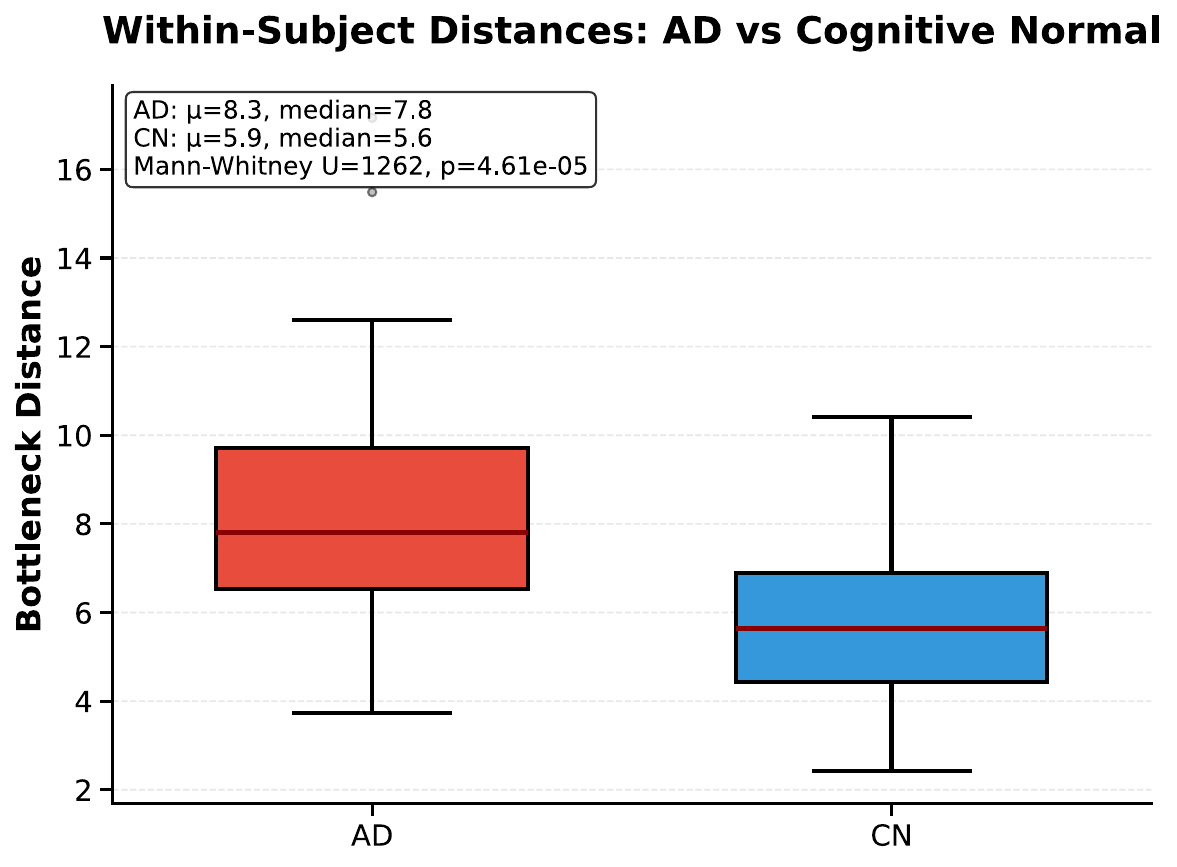}
\\
\end{tabular}
\caption{\label{fig:within_between_stability}Longitudinal stability analysis of \(H_2\) persistence. Boxplots show within-subject (\(n=82\)) versus subject-level between-subject distances, where \(B_i\) is the median baseline-baseline distance from subject \(i\) to all other subjects. Statistical inference uses paired Wilcoxon (\(B_i > W_i\), greater): \(W=3403\), \(p = 1.83 \times 10^{-15}\). Within-subject median \(d_B = 6.66\) mm\(^2\); subject-level between-subject median \(\tilde{B} = 26.78\) mm\(^2\); median per-subject ratio \(R_i = 4.32\times\) (median-based \(B_i\)).}
\end{figure*}

To avoid dependence among the \(13{,}284\) pairwise distances, we conducted statistical inference at the subject level. For each subject \(i\), we computed the within-subject distance \(W_i\) (baseline vs follow-up) and a between-subject summary \(B_i\) using two methods: (1) \textbf{median} and (2) \textbf{mean} of baseline-baseline distances from subject \(i\) to all other subjects. Element-wise ratios \(R_i = B_i / W_i\) were tested using paired Wilcoxon signed-rank tests (\(B_i > W_i\), one-sided) (Table \ref{tab:stability}).

\begin{table*}[htbp]
\caption{\label{tab:stability}Within-subject vs between-subject bottleneck distances for longitudinal stability analysis}
\centering
\small
\begin{tabular}{lcccccc}
\toprule
Comparison & N subjects & Median \(d_B\) (mm\(^2\)) & Mean \(d_B\) (mm\(^2\)) & Median \(R_i\) & Mean \(R_i\) & Wilcoxon\\
\midrule
(\(W_i\)) & 82 & 6.66 & 7.09 & - & - & -\\
median-based \(B_i\) & 82 & 26.78 & 32.83 & 4.32\(\times\) & 4.91\(\times\) & \(W=3403\)\\
mean-based \(B_i\) & 82 & 31.02 & 35.81 & 5.12\(\times\) & 5.55\(\times\) & \(W=3403\)\\
\bottomrule
\end{tabular}
\end{table*}

Every subject's follow-up scan was closer to its own baseline than to any other subject's baseline: both summary methods achieved the maximum possible Wilcoxon statistic (\(W=3403 = \sum_{i=1}^{82} i\), \(p = 1.83 \times 10^{-15}\)), with median per-subject ratios of \(4.32\times\) (median-based) and \(5.12\times\) (mean-based). This shows that the \(H_2\) diagrams retain subject specificity over six-month intervals despite disease-related change.

Stratification by diagnostic group (Tables \ref{tab:within_by_group} and \ref{tab:ad_cn_stats}) showed a \(41\%\) increase in within-subject variability for AD subjects (\(8.32\) vs \(5.91\) mm\(^2\)). This aligns with accelerated and heterogeneous CSF ex-vacuo expansion in disease progression.

\begin{table}[htbp]
\caption{\label{tab:within_by_group}Within-subject bottleneck distances by diagnostic group}
\centering
\small
\begin{tabular}{lccc}
\toprule
Group & N subjects & Median \(d_B\) & Mean \(d_B\)\\
\midrule
CN & 42 & 5.63 & 5.91\\
AD & 40 & 7.80 & 8.32\\
\bottomrule
\end{tabular}
\end{table}

\begin{table}[htbp]
\caption{\label{tab:ad_cn_stats}Within-subject variability differs between AD and CN.}
\centering
\small
\begin{tabular}{lccc}
\toprule
Comparison & Test statistic & \(p\)-value & Effect size\\
\midrule
AD vs CN & Mann-Whitney \(U=1262\) & \(4.61 \times 10^{-5}\) & Cohen's \(d=0.95\)\\
\bottomrule
\end{tabular}
\end{table}
\subsection*{Topological Features Outperform CSF Volumetrics}
\label{sec:org34e271b}
Cross-sectionally, CSF\% provided a simple volumetric baseline. Using the same classifier (PCA95+sLDA), CSF\% achieved ROC-AUC = 0.840 and balanced accuracy = 0.740 in CN vs AD, compared with ROC-AUC = 0.895 and balanced accuracy = 0.814 for Pipeline 1.

Longitudinally, relative CSF change performed near chance in this cohort (ROC-AUC = 0.569, \(p = 0.28\), Cohen's \(d = 0.19\); see Supplement), whereas \(H_2\) bottleneck distances performed better (ROC-AUC = 0.751, Cohen's \(d = 0.95\)). Together, these comparisons show that topology captures geometric information beyond scalar volumetric measures.
\section{Discussion}
\label{sec:org12c2fe6}
Structural MRI morphometry is usually built around nonlinear warping to a common template, atlas-based regional summaries, or mass-univariate voxel-wise inference \citep{Ashburner2000VBM}. These frameworks have been highly influential, but they also impose a particular assumption: that brain change is best measured after standardizing anatomy into a shared coordinate system or reducing it to scalar regional volumes. Our results challenge that assumption. By treating homological cycles as the primitive objects of analysis, we show that biologically meaningful atrophy signals can be measured directly in subject space from standard tissue masks, without nonlinear registration, while still supporting localization, classification, and longitudinal tracking.
\subsubsection*{Homological Cycles as Primitives.}
\label{sec:orgc663bb2}
The central methodological shift in this work is the adoption of homological cycles (loops and voids), rather than voxels, as primitive objects for feature engineering, and topological summaries (landscape \(L^{1}\) norm and bottleneck distance) as the basis for statistical inference. Within this framework, we show that brain atrophy carries a topological signature, i.e., a reshaping of spatial arrangement that is measurable, interpretable, and diagnostically informative from structural MRI alone.

This builds on a growing body of work adopting persistent homology for medical imaging \citep{Singh2023,Singh_2024,Salch_2021}. In AD-related neurodegeneration, prior approaches demonstrated separation from controls using PH cortical thickness maps \citep{Pachauri2011} and coordinate-free descriptors from hippocampal shape \citep{Hofer2017}. However, these required computationally intensive surface reconstructions, were limited to small samples, or achieved modest accuracy \citep{SaadatYazdi2021}. A PH framework that works without template normalization, starting from standard clinical preprocessing of structural images, was missing.

The present work fills this gap. The choice of PH filtration plays a key role in the design of biomarkers, uncovering brain signatures of normal and disease-related atrophy. The filtration determines what aspect of brain geometry the cycles encode: EDT superlevel filtrations (Pipeline 1) produce persistent \(H_1\) cycles whose lifetimes reflect tissue depth and thinning, while \(\alpha\)-complex filtrations (Pipeline 2) produce persistent \(H_2\) cycles whose lifetimes reflect cavity size and CSF expansion.

These correspondences directly link homology to the biophysics of brain atrophy and its morphological evolution in time. This distinguishes the present approach from prior TDA work that used intensity-based or anatomically ungrounded filtrations \citep{Bruningk2021}.
\subsubsection*{Anatomical Correspondence Without Warping.}
\label{sec:org45620d2}
A dominant challenge in structural neuroimaging is the need to design features from images that can be compared. All brains share the same gross anatomy, but their local geometry---sulcal folding, gyri patterns and ventricle shapes---varies considerably.

In VBM, nonlinear registration (warping) to standard templates and kernel smoothing are the de facto solutions to this problem. This introduces well-documented biases: registration accuracy varies across brain regions and disease severity, atrophied brains are harder to register, and the choice of template itself influences results \citep{Diaz_de_Grenu_2013,Zhou_2022,Giuliani_2005}. After warping, modulation\footnote{Modulation by the Jacobian determinant of the deformation.} and smoothing, the interpretation of ``tissue densities'' can still be challenging. On the other hand, ROI-based volumetry avoids some of these problems but requires atlas-based parcellation and reduces rich spatial information to single scalar volumes \citep{Yaakub_2020}.

Our pipelines take a different route. They operate in subject space without requiring nonlinear registration.

Pipeline 1 integrates local measurements into global summary statistics (\(L^{1}\) curves). The choice to analyze 2D slices rather than full 3D volumes is deliberate: it forces locality. Each slice provides a local readout of the tissue's geometric organization at a known position along an anatomical axis. Furthermore, because the brain is rigidly aligned to a standard orientation, normalized slice positions correspond to approximately the same anatomical territory across subjects. The aggregation of these slice-wise measurements into a single curve yields a global morphometric profile that is effectively anchored by the gross structural anatomy shared at each normalized position.

For instance, it is visible in the mean \(L^1\) curves (Fig. \ref{fig:l1_curves_all_groups}): at the midsagittal plane, the \(L^1\) norm is consistently low---the dominant midline structures (corpus callosum, brainstem and cerebellar vermis) are thick and largely solid, while the enclosed midline CSF spaces are narrow and bounded by thin tissue, yielding only low-persistence \(H_1\) cycles. Conversely, \(L^1\) peaks mark slices rich in persistent cycles, where tissue surrounds the ventricular system and deep sulci.

The slice-wise decomposition thus reduces the comparison problem from matching two 3D objects with wildly different local geometry to comparing 1D curves that respect both the global similarity and local uniqueness of individual brains.

Pipeline 2 uses \(\alpha\)-complex filtrations to summarize the geometry of CSF-bounded cavities from a point-cloud representation of the CSF complement \citep{Edelsbrunner2009AlphaS}. As the filtration parameter increases, the alpha shape gradually fills the space occupied by the CSF system. In this setting, persistent \(H_2\) cycles correspond to enclosed cavities bounded by ventricles, sulci, and cisterns. Their persistence therefore reflects the size and distribution of CSF spaces that enlarge as tissue is lost. Comparing these \(H_2\) persistence diagrams across visits with the bottleneck distance yields a direct measure of longitudinal morphological change. Because these measurements are coordinate-free and isometry-invariant, they do not depend on parcellations or template warping.
\subsubsection*{Contributions to AD-related Morphometry.}
\label{sec:org801d8b3}
Our cross-sectional results demonstrate that topological features capture AD-related atrophy beyond what is explained by age or brain size. Stratification by age and total intracranial volume preserved CN vs AD separation across all demographic strata (Cohen's \(d\) range: \(0.51\)--\(1.51\)), confirming that the signal reflects disease pathology rather than normal aging or head-size variation. Moreover, the \(L^1\) curves provide a measure that both discriminates disease groups and localizes atrophy-related effects.
Unlike deep learning approaches that achieve comparable accuracy at the cost of diluting interpretability and reproducibility \citep{Wen2020,Ebrahimighahnavieh2020}, our features maintain direct correspondence to anatomical structures.

The localization of peak effects to medial temporal structures (hippocampus, parahippocampal gyrus, amygdala) aligns with known vulnerable regions of AD pathology. These regions show the early neurofibrillary tangle accumulation in Braak staging and exhibit selective vulnerability \citep{Braak_1995,Dickerson_2008,Stoub_2014}. The fact that our data-driven localization independently recovers these anatomically plausible regions supports the biological validity of the topological features.

Our classification performance (ROC-AUC = \(0.895\pm 0.034\); balanced acc. = \(0.814 \pm 0.048\); CN vs AD) compares favorably to prior PH-based approaches. Pachauri et al. achieved ROC-AUC = \(0.80\) using persistent homology of cortical thickness maps \citep{Pachauri2011}. Saadat-Yazdi et al. reported balanced acc. of \(0.76 \pm 0.02\) from gray, white matter and surface point cloud topological features \citep{SaadatYazdi2021}. Our framework achieves stronger discrimination while requiring only standard tissue segmentation, no surface reconstruction, and no nonlinear registration. These results were obtained from single-modality T1w structural MRI without additional imaging biomarkers, showing the information content available from brain morphometry alone.

The longitudinal analysis shows that these topological descriptors retain subject specificity over time. Each follow-up scan remained closer to its own baseline than to any other subject's baseline, with between-subject distances \(5.12\) times greater than within-subject distances. This strong within-subject stability makes the bottleneck distance sensitive to short-interval disease-related change rather than merely to inter-individual anatomical variation.

AD subjects also showed \(41\%\) greater within-subject variability than CN subjects over six-month intervals (Cohen's \(d = 0.95\)). This pattern is consistent with faster and more heterogeneous CSF expansion during disease progression.

The comparison with CSF\% provides the simplest scalar baseline from the same tissue masks. Cross-sectionally, Pipeline 1 outperformed CSF\% (ROC-AUC = \(0.895 \pm 0.034\) versus \(0.840 \pm 0.039\)). Longitudinally, relative CSF change performed near chance over six-month intervals (ROC-AUC = \(0.569\)), whereas the \(H_2\) bottleneck distance remained discriminative (ROC-AUC = \(0.751\)). This suggests that, at short timescales, total CSF volume changes are small relative to measurement noise, while the spatial redistribution of CSF spaces carries usable disease information.
\subsubsection*{Practical Considerations.}
\label{sec:org1c3f77d}
The entry point for both pipelines is a standard tissue segmentation mask. Both pipelines run on commercial hardware at the scales presented in this study (\(n=929\) scans after quality control). This is relevant for clinical settings, where computational burden and preprocessing complexity remain barriers to adoption.

In both pipelines, the effect of segmentation error is bounded. The stability theorem of persistent homology \citep{CohenSteiner2007Stability,Chazal_2013} ensures that a boundary perturbation of at most \(\varepsilon\) shifts the persistence diagram by at most \(\varepsilon\) in bottleneck distance. All segmentations also underwent systematic visual quality control, and scans with gross tissue misclassification, skull-stripping failures, or discontinuities in the CSF mask were excluded (see Supplement).

We emphasize that the pipelines target atrophy itself, not Alzheimer's disease specifically. We chose AD as testing ground because it is a well-characterized neurodegenerative disease with pronounced atrophy. The synthetic validation supports this generality: with controlled erosion experiments, both pipelines responded strongly to the induced volume loss (Spearman \(\rho = -0.76\) for Pipeline 1; \(\rho = 0.928\) for Pipeline 2). Because the features are tied to general geometric consequences of tissue loss---not to AD-specific molecular pathology---the framework should extend to any condition involving parenchymal volume loss and CSF expansion.
\subsubsection*{Limitations and Future Directions.}
\label{sec:orga98b29b}
Persistent homology does not identify a unique geometric representative for each cycle. This means that individual persistence pairs should not be interpreted as anatomically unique loops or cavities. However, this abstraction is also part of the strength of the method: it focuses on topological invariants of the space rather than on voxel-level geometry, making the resulting descriptors less sensitive to local perturbations. In our setting, this trade-off is appropriate because the quantities used for inference are aggregate summaries (landscape \(L^1\) norms and bottleneck distances), which are designed to capture distributional structure rather than the exact geometry of any single cycle. For Pipeline 2, this does not eliminate correspondence across scans: the bottleneck matching still links prominent \(H_2\) features across visits, especially for large stable CSF-bounded cavities. However, that correspondence is topological rather than anatomically unique, so the diagrams should be read as tracking the geometry of cavity systems rather than assigning one fixed label to each individual sulcus or ventricle.

Pipeline 2's computational burden scales with point cloud size. At \(1\) mm\(^3\) resolution, \(\alpha\)-complex construction operates on \(\sim\!3 \times 10^{6}\) points per subject, requiring \(\sim\!140\)s per scan. Resampling to \(2\) mm\(^3\) reduces the point cloud to \(\sim\!4 \times 10^{5}\) points, cutting computation to \(\sim\!16\)s (\(8.7\times\) speedup). This ablation fully preserved within-subject stability: every subject's follow-up remained closer to its own baseline than to any other subject (Wilcoxon \(W = 3403\), \(p = 1.83 \times 10^{-15}\)), with a median per-subject ratio of \(4.09\times\) (vs \(5.12\times\) at \(1\) mm\(^3\)). AD vs CN discrimination remained significant (ROC-AUC = \(0.721\)), though with a reduced effect size (Cohen's \(d = 0.59\) vs \(0.95\) at \(1\) mm\(^3\)). The \(2\) mm\(^3\) resolution thus offers a practical trade-off for large-scale studies where longitudinal fingerprinting is the primary goal, while \(1\) mm\(^3\) resolution maximizes sensitivity for group-level discrimination.

The only source of stochasticity in our framework is the rigid registration used in Pipeline 1. We fixed the random seed and enforced single-threaded execution to make this step reproducible. Even without these controls, the residual registration-induced variation had little effect on final classification performance. Pipeline 2 is fully deterministic.

The topological features presented here are structural by design. The modest CN vs MCI discrimination (ROC-AUC = \(0.699\)) reflects the limited macroscopic atrophy in early disease stages, highlighting a ceiling of single-modality structural imaging for early detection. At the MCI stage, disease-related change is still dominated by molecular and functional events: amyloid and tau accumulation, synaptic dysfunction, and early network reorganization. Parenchymal loss has barely begun to leave a macroscopic footprint on T1-weighted MRI. The limited information available from T1w alone at this stage therefore bounds what any structural morphometric method can recover, topological or otherwise.

A natural next step is to combine PH-derived geometric features with modalities that probe earlier stages of disease, including amyloid or tau PET and functional signals from fMRI, fNIRS, or EEG. Persistence diagrams and their summaries can be integrated with such measurements through standard multimodal fusion strategies. More broadly, because the pipelines target general consequences of tissue loss rather than AD-specific anatomy, they should also transfer to other disorders with progressive atrophy, including traumatic brain injury, Parkinson's disease, and multiple sclerosis.
\subsubsection*{Conclusions.}
\label{sec:org8e70ab5}
This paper argues that brain atrophy is not just a loss of quantity but a transformation of shape. Persistent homology reveals a level of structural description between raw geometry and scalar reduction---the level of topology---that is both mathematically rigorous and biologically informative.

Homology-based morphometry offers an alternative to existing morphometric frameworks that can target different aspects of brain geometry through different filtrations while maintaining low stochasticity. By treating homological cycles as primitive objects, our approach minimizes preprocessing while preserving interpretability.

We applied this framework to both cross-sectional and longitudinal study designs, achieving strong single-modality classification and sensitive short-term longitudinal tracking from standard T1w structural MRI.

\bibliography{paper1--topo-adni}

\newpage

\onecolumn
\section{Supplementary Material}
\label{sec:orgb4fec8b}

In this section we provide complementary analyses for our pipelines. We discuss correlations of our features with clinical measures to better understand the possible clinical use cases of our framework. We then provide a stratification analysis for Pipeline 1; this is needed to assess confounding effects with well-known covariates like brain size and age. Then we compare our work to CSF\%, the simplest volumetric measurement that can be computed from the same inputs as our pipelines and a strong baseline for Alzheimer's disease classification. We then show the full Pipeline 1 analysis results including MCI subtypes (EMCI and LMCI). Finally we discuss the quality control at the preprocessing stage.
\subsection{Correlation with Clinical Measures}
\label{sec:orgd40ca16}

To establish clinical relevance, we correlated \(L^1\) curves (Pipeline 1) with cognitive and volumetric measures from ADNIMERGE. We matched clinical assessments to MRI sessions within a \(\pm\)90-day window to balance temporal proximity and sample size. Overall, \(92.7\%\) of subjects were matched (\(709\) of \(765\)), and \(85.9\%\) of matches occurred within \(\pm\)30 days (mean temporal distance of 16.4 days). This window accommodates visit scheduling variation while limiting clinical change between assessments.

\begin{figure}[H]
\centering
\includegraphics[width=0.6\textwidth]{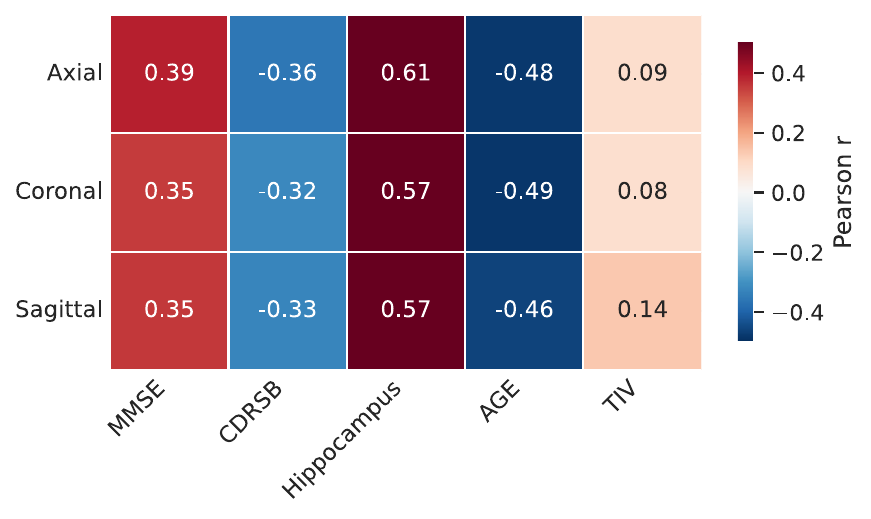}
\caption{\label{fig:correlation_matrix}Pearson correlations between \(L^1\) curves AUC and clinical measures. Sample sizes vary due to missing clinical data: Age (\(n=709\)), TIV (\(n=704\)), MMSE (\(n=700\)), CDR-SB (Clinical Dementia Rating--Sum of Boxes; \(n=699\)), Hippocampus (\(n=641\)).}
\end{figure}

Topological features correlate with established AD biomarkers (Fig. \ref{fig:correlation_matrix}). The strongest correlation is with hippocampal volume (\(r = 0.61\), \(p = 5.98 \times 10^{-67}\)).

Cognitive measures show moderate correlations: MMSE (\(r = 0.39\), \(p = 1.32 \times 10^{-26}\)) and CDR-SB (\(r = -0.36\), \(p = 9.28 \times 10^{-23}\)), in the expected direction (higher \(L^1\)-curve AUC correlates with better cognition). Age correlates negatively (\(r = -0.48\), \(p = 1.86 \times 10^{-42}\)). TIV correlation is weak (\(r = 0.14\), \(p = 2.80 \times 10^{-4}\)), suggesting minimal head-size confounding. Taken together, these associations indicate that the topological features track clinically relevant atrophy signal rather than merely reflecting age or global brain size.

\begin{figure*}
\centering
\includegraphics[width=0.95\textwidth]{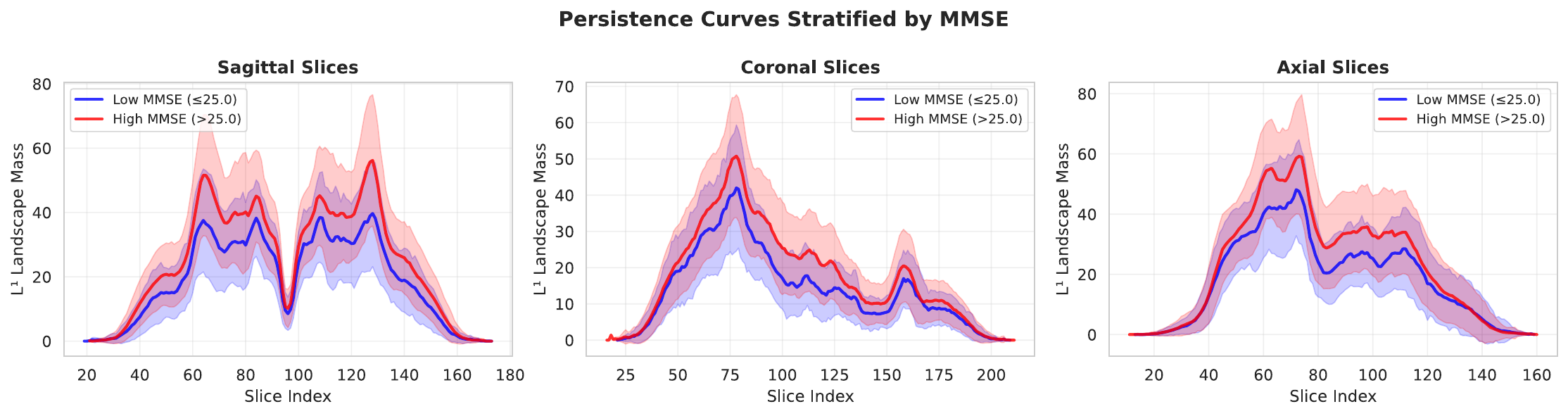}
\caption{\label{fig:curves_mmse}\(L^1\) curves (mean \(\pm\) SD) stratified by MMSE. High MMSE subjects (\(>25\), red) show higher \(L^1\) landscape AUC across all anatomical views compared to low MMSE subjects (\(\le25\), blue).}
\end{figure*}

To visualize the relationship between topological features and cognitive status, we stratified subjects by MMSE score and computed mean \(L^{1}\) curves for high (\(>25\)) versus low (\(\le 25\)) groups (Fig. \ref{fig:curves_mmse}).

\subsection{Stratification Analysis}
\label{sec:org0924e4b}

To assess whether group differences in topological features are confounded by age or total intracranial volume (TIV), we performed stratified analyses controlling for both covariates. Age correlates with brain atrophy in both healthy aging and AD, while TIV represents brain size that could influence absolute feature magnitudes.
\subsubsection{Methodology}
\label{sec:org62f844c}

We stratified the 3T cross-sectional cohort (CN and AD subjects, \(n=334\); 201 CN, 133 AD) using median splits to create two balanced bins for each covariate:

\begin{table}[htbp]
\caption{\label{tab:strata}Strata}
\centering
\small
\begin{tabular}{c l l}
\toprule
Variable & Category & Definition\\
\midrule
Age (years) & Young & \(< 75.0\) years\\
 & Old & \(\geq 75.0\) years\\
\midrule
Total intracranial & Low TIV & \(< \SI{1490000}{mm^3}\)\\
volume (TIV) & High TIV & \(\geq \SI{1490000}{mm^3}\)\\
\bottomrule
\end{tabular}
\end{table}

This 2\(\times\)2 stratification produces four subgroups for analysis. Within each stratum, we examined:
\begin{enumerate}
\item Mean \(L^1\) curves (\(\pm\) SD) across anatomical axes to assess whether disease-related morphological patterns persist after controlling for confounds.
\item Effect sizes (Cohen's \(d\)) on total \(L^1\)-curve AUC (summed across axes) to quantify the magnitude of group differences within each stratum.
\end{enumerate}
\subsubsection{Results}
\label{sec:org18590ee}

\begin{figure*}
\centering
\includegraphics[width=0.95\textwidth]{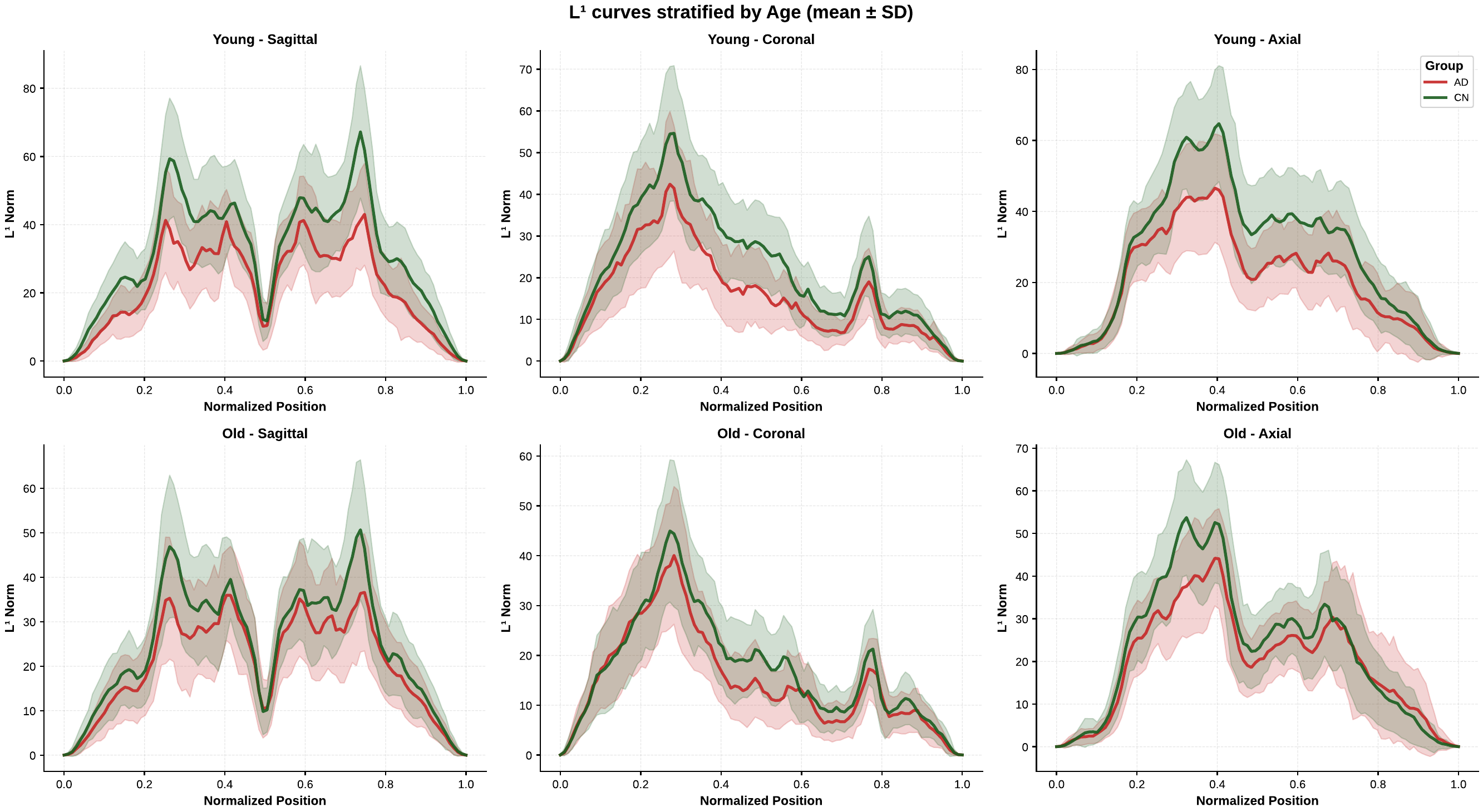}
\caption{\label{fig:l1_age_stratified}\(L^1\) curves stratified by age. Mean \(\pm\) SD for CN (green) and AD (red) across sagittal, coronal, and axial views. Top: Young (\(<\)75 years). Bottom: Old (\(\geq\)75 years). CN-AD separation persists in both strata.}
\end{figure*}

\begin{figure*}
\centering
\includegraphics[width=0.95\textwidth]{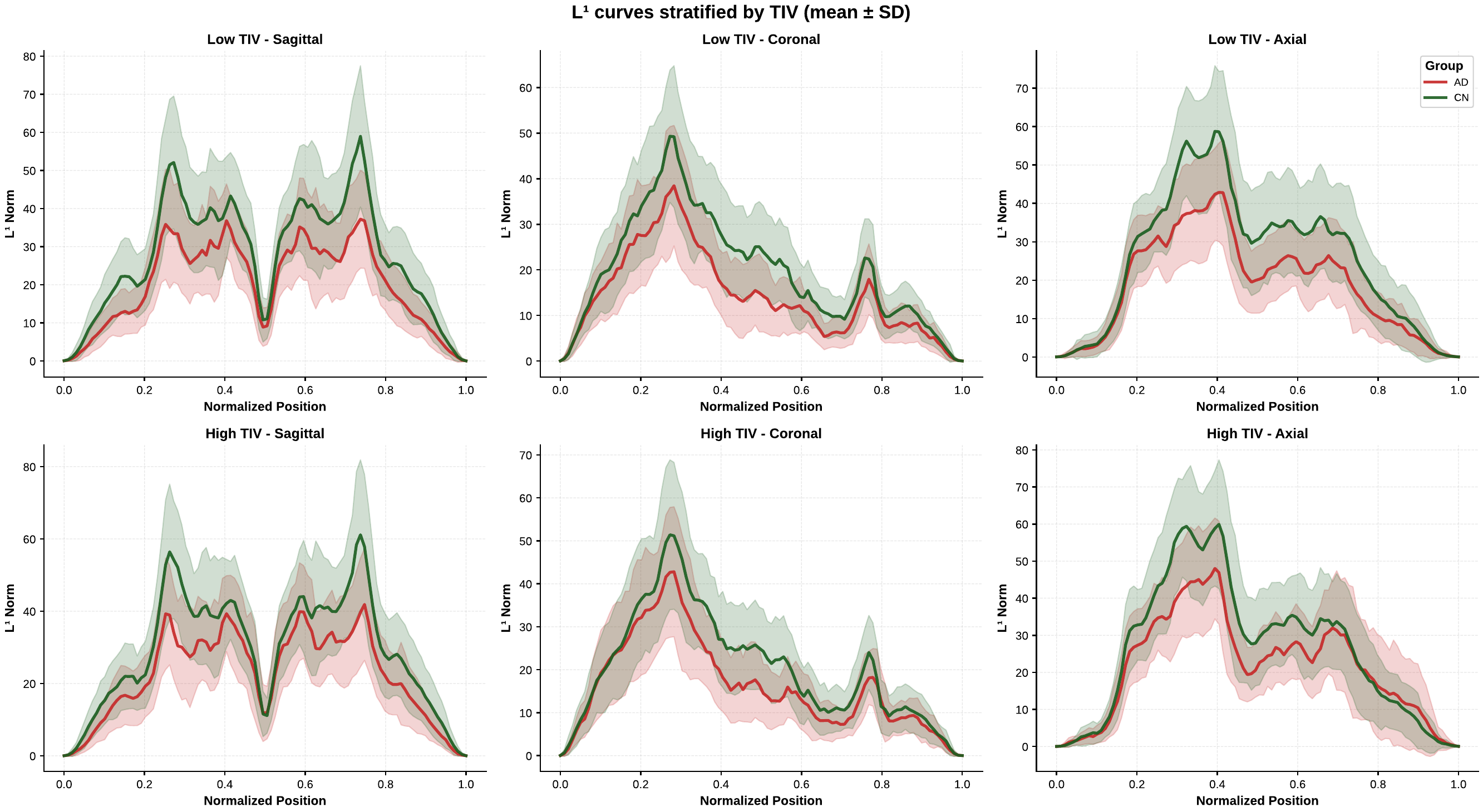}
\caption{\label{fig:l1_tiv_stratified}\(L^1\) curves stratified by TIV. Mean \(\pm\) SD for CN (green) and AD (red) across sagittal, coronal, and axial views. Top: Low TIV. Bottom: High TIV. CN-AD separation persists in both strata.}
\end{figure*}

Differences in \(L^1\) curves between CN and AD subjects persisted after stratification by both age and TIV (Figures \ref{fig:l1_age_stratified}, \ref{fig:l1_tiv_stratified}). Within every stratum, AD subjects showed lower \(L^1\)-curve AUC than CN across all three anatomical views, suggesting these group differences are not driven by normal aging or brain size.

Cohen's \(d\) effect sizes on total \(L^1\)-curve AUC (summed across axes) were computed for each Age \(\times\) TIV stratum using independent samples \(t\)-tests (Table \ref{tab:effect_sizes}). Effect sizes ranged from \(d = 0.51\) (Old/High TIV, \(p = 1.96 \times 10^{-2}\)) to \(d = 1.51\) (Young/Low TIV, \(p = 7.11 \times 10^{-10}\))

\begin{table*}[htbp]
\caption{\label{tab:effect_sizes}Effect sizes for CN vs AD within Age \(\times\) TIV strata. Total AUC = area under \(L^1\) curves summed across all axes.}
\centering
\small
\begin{tabular}{lrrccc}
\toprule
Stratum & \(N_{CN}\) & \(N_{AD}\) & AUC diff. (\%) & Cohen's \(d\) & \(p\)-value\\
\midrule
Young, Low TIV & 60 & 32 & -41.3 & 1.51 & \(7.11 \times 10^{-10}\)\\
Young, High TIV & 54 & 25 & -34.4 & 1.46 & \(5.24 \times 10^{-8}\)\\
Old, Low TIV & 41 & 34 & -25.7 & 1.17 & \(3.39 \times 10^{-6}\)\\
Old, High TIV & 46 & 42 & -10.9 & 0.51 & \(1.96 \times 10^{-2}\)\\
\bottomrule
\end{tabular}
\end{table*}

In conclusion, this analysis shows that CN vs AD separation persists across age and TIV strata.
\subsection{CSF\% Baseline}
\label{sec:orga82efa7}

To assess whether topological features capture information beyond volumetrics, we computed CSF\% (CSF/(CSF + GM + WM)) as a baseline. Masks were resampled to 1mm\(^3\) isotropic voxels before counting.
\subsubsection{Cross-Sectional Analysis}
\label{sec:orgbad5487}
CSF\% discriminated CN from AD: CN \(24.4\% \pm 2.6\%\) versus AD \(28.2\% \pm 2.6\%\) (Cohen's \(d = 1.45\), \(p = 2.01 \times 10^{-31}\)). To ensure a fair comparison, we evaluated CSF\% using the same classifiers as Pipeline 1 (PCA95+sLDA, RBF-SVM, XGBoost; 5-fold stratified CV). The best-performing model (PCA95+sLDA) achieved ROC-AUC = \(0.840 \pm 0.039\) and balanced accuracy = \(0.740 \pm 0.037\).

\begin{figure*}
\centering
\includegraphics[width=0.5\textwidth]{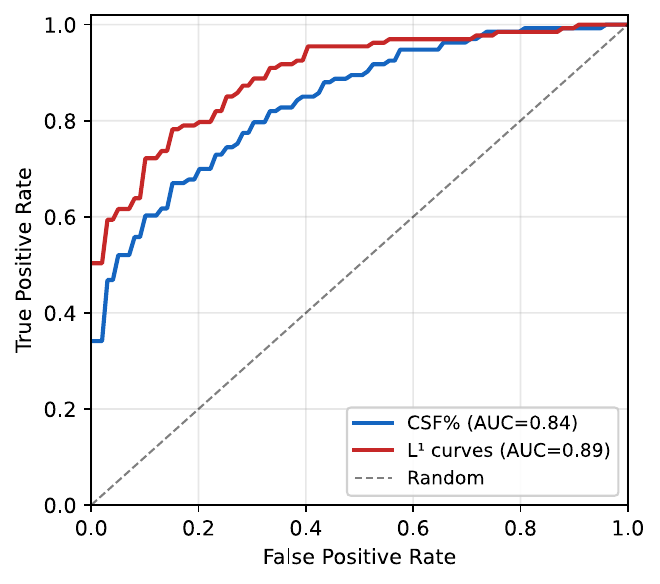}
\caption{\label{fig:roc-plot}ROC comparison between CSF\% and \(L^1\) curves.}
\end{figure*}

CSF\% achieved lower classification performance than Pipeline 1 (Table \ref{tab:csf_vs_topo}, Fig. \ref{fig:roc-plot}). Within groups, CSF\% and \(L^1\) AUC correlate strongly (CN: \(r = -0.79\); AD: \(r = -0.75\)), but topology captures spatial organization that scalar ratios miss.

\begin{table}[htbp]
\caption{\label{tab:csf_vs_topo}Cross-sectional classification: CSF\% vs topological features (PCA95+sLDA for both; 5-fold CV)}
\centering
\small
\begin{tabular}{lcc}
\toprule
Metric & CSF\% & Pipeline 1\\
\midrule
ROC-AUC & \(0.840 \pm 0.039\) & \(0.895 \pm 0.034\)\\
Balanced Acc & \(0.740 \pm 0.037\) & \(0.814 \pm 0.048\)\\
\bottomrule
\end{tabular}
\end{table}
\subsubsection{Longitudinal Analysis}
\label{sec:orge3fa03b}

For the 1.5T longitudinal cohort (\(n=82\)), relative CSF change (\(\Delta\text{CSF\%} = \text{CSF\%}_{t_2} - \text{CSF\%}_{t_1}\)) did not discriminate AD from CN progression:

\begin{table}[htbp]
\caption{\label{tab:csf_longitudinal_comparison}Longitudinal discrimination: relative CSF change vs \(H_2\) bottleneck distance}
\centering
\small
\begin{tabular}{l c c}
\toprule
Metric & ROC-AUC & \(p\)-value\\
\midrule
Relative CSF change & 0.569 & 0.28\\
\(H_2\) bottleneck & 0.751 & < 0.0001\\
\bottomrule
\end{tabular}
\end{table}

Within-subject CSF changes were small relative to between-subject variability (\(p = 0.28\)). At this timescale, CSF\% change is dominated by measurement noise. \(H_2\) persistence, by contrast, captures geometric changes in CSF cavities beyond total volume.
\subsection{MCI Subtypes}
\label{sec:orgd288081}
The original ADNI1 protocol classified subjects with intermediate cognitive decline as a single MCI group. Beginning with ADNI-GO, this category was subdivided into Early MCI (EMCI) and Late MCI (LMCI) based on education-adjusted delayed paragraph recall scores (Logical Memory II, Wechsler Memory Scale--Revised), with EMCI capturing milder impairment (approximately \(1.0\) SD below norms) and LMCI corresponding to the more severe threshold used in ADNI1 (approximately \(1.5\) SD below norms) \citep{Aisen_2010}. Because our 3T dataset spans ADNI1, ADNI-GO and ADNI2, it contains subjects labeled under both schemes: MCI (ADNI1, \(n=97\)) alongside EMCI (\(n=216\)) and LMCI (\(n=118\)) from ADNI-GO/ADNI2. LMCI subjects meet the same memory impairment threshold as ADNI1 MCI; the main analysis uses only the ADNI1 MCI cohort to avoid cross-phase confounds, but the comparable LMCI results below support generalizability. Here we report the EMCI and LMCI subtype analyses.

\begin{table*}[htbp]
\caption{\label{tab:demographics_subgroups}Demographics by diagnostic group including MCI subtypes (3T dataset).}
\centering
\small
\begin{tabular}{r c c c c}
\toprule
\textbf{Group} & \textbf{N (count)} & \textbf{Age (mean\(\pm\)SD)} & \textbf{MMSE (mean\(\pm\)SD)} & \textbf{TIV (mean\(\pm\)SD)}\\
\midrule
CN & 201 & 74.48\(\pm\) 6.70 & 28.94\(\pm\) 1.32 & \num{1424096}\(\pm\)\num{133255}\\
EMCI & 216 & 70.74\(\pm\) 7.04 & 28.42\(\pm\) 1.55 & \num{1435360}\(\pm\)\num{132797}\\
LMCI & 118 & 72.03\(\pm\) 8.07 & 27.50\(\pm\) 1.94 & \num{1433877}\(\pm\)\num{140014}\\
MCI & 97 & 75.22\(\pm\) 7.92 & 27.12\(\pm\) 2.31 & \num{1462362}\(\pm\)\num{129617}\\
AD & 133 & 75.98\(\pm\) 7.88 & 20.14\(\pm\) 4.86 & \num{1428998}\(\pm\)\num{157331}\\
\bottomrule
\end{tabular}
\end{table*}

Topological features show modest discrimination between MCI subtypes. We applied the same classification framework used for CN vs AD to MCI subtype comparisons, selecting the best-performing model (XGBoost) based on ROC-AUC.

\begin{figure*}
\centering
\includegraphics[width=0.95\textwidth]{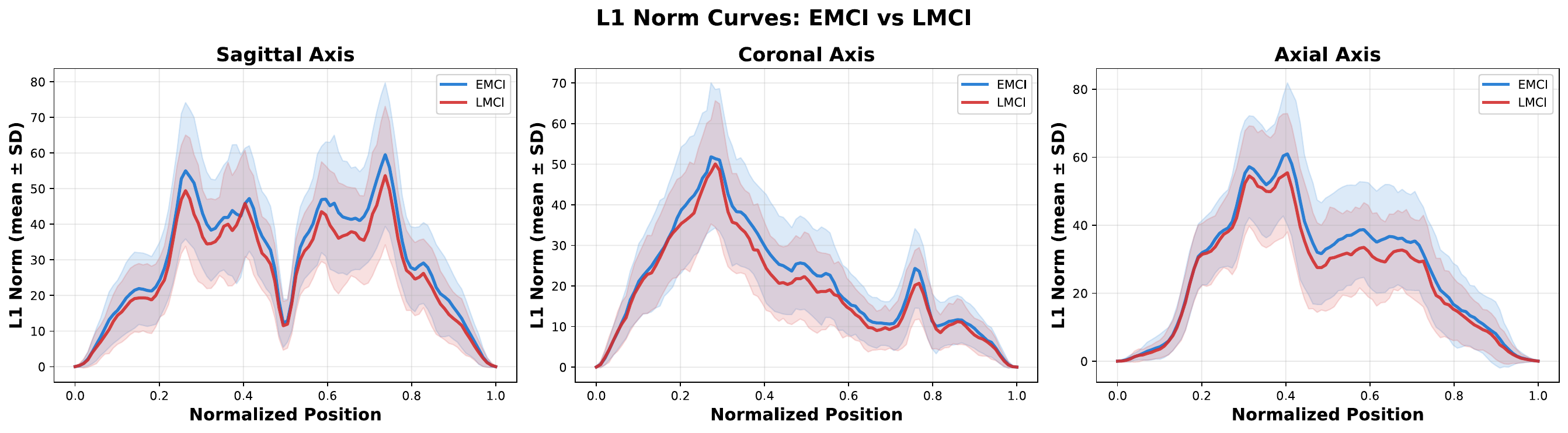}
\caption{\label{fig:emci_lmci_comparison}\(L^1\) curves comparing EMCI (blue) versus LMCI (red) across anatomical orientations. EMCI subjects show higher \(L^1\)-curve values, reflecting less advanced tissue degradation. Shaded regions show \(\pm\) SD.}
\end{figure*}

EMCI subjects maintain higher \(L^1\)-curve values than LMCI across all anatomical views, particularly in sagittal and axial orientations (Fig. \ref{fig:emci_lmci_comparison}).

\begin{table*}[htbp]
\caption{\label{tab:mci_classification}Classification performance for MCI subtype comparisons (XGBoost, 5-fold CV)}
\centering
\small
\begin{tabular}{l c c c c}
\toprule
Comparison & N & Balanced Accuracy & F1 Score & ROC-AUC\\
\midrule
EMCI vs LMCI & 334 & 0.571 \(\pm\) 0.048 & 0.368 \(\pm\) 0.085 & 0.604 \(\pm\) 0.051\\
CN vs LMCI & 319 & 0.635 \(\pm\) 0.019 & 0.497 \(\pm\) 0.046 & 0.683 \(\pm\) 0.018\\
LMCI vs AD & 251 & 0.688 \(\pm\) 0.042 & 0.709 \(\pm\) 0.035 & 0.757 \(\pm\) 0.046\\
EMCI vs AD & 349 & 0.797 \(\pm\) 0.056 & 0.748 \(\pm\) 0.073 & 0.884 \(\pm\) 0.044\\
\bottomrule
\end{tabular}
\end{table*}

Classification across MCI subtypes (Table \ref{tab:mci_classification}) shows strongest discrimination between adjacent disease stages in LMCI vs AD (ROC-AUC = \(0.757 \pm 0.046\)) and modest discrimination between EMCI and LMCI (ROC-AUC = \(0.604 \pm 0.051\)), in line with the clinical similarity of disease stages.
\subsection{Quality Control}
\label{sec:orge38ebe7}

\begin{figure*}
\centering
\includegraphics[width=0.9\textwidth]{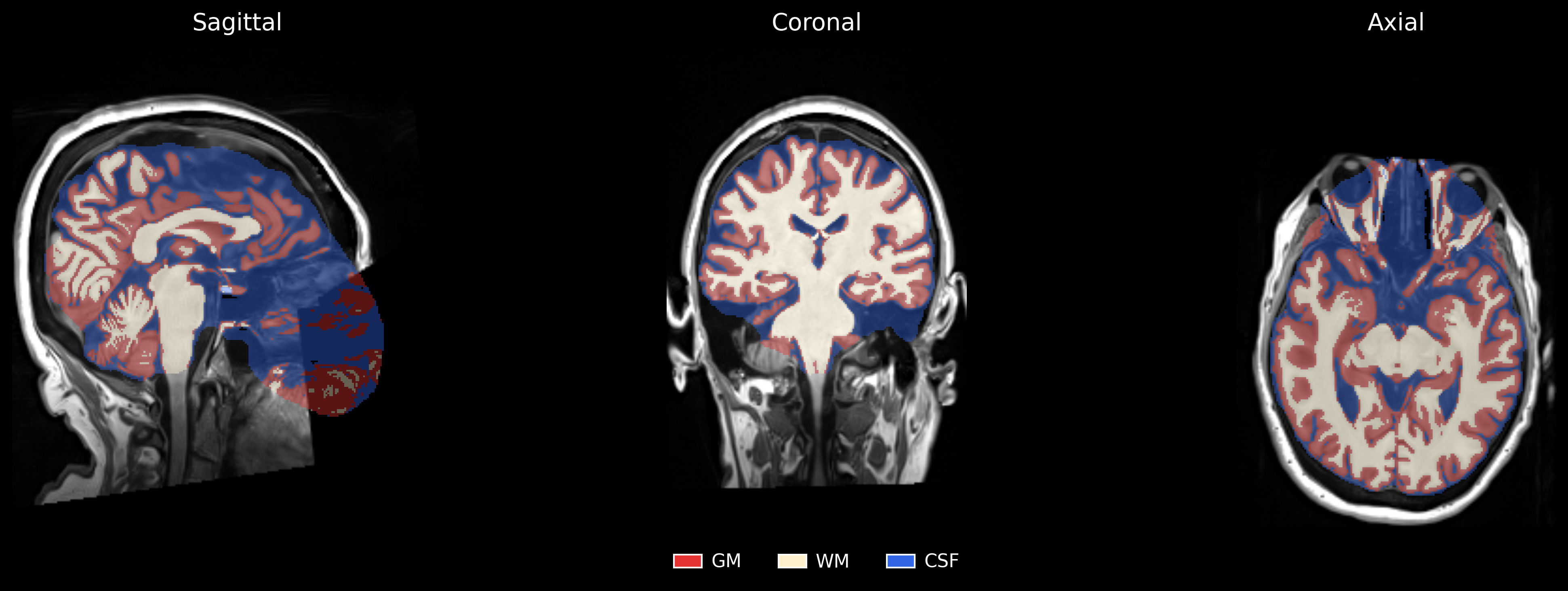}
\caption{\label{fig:qc_gross_failure}Excluded segmentation showing tissue classes overlaid on defaced T1w: GM (red), WM (cream), CSF (blue).}
\end{figure*}

\begin{figure*}
\centering
\includegraphics[width=0.9\textwidth]{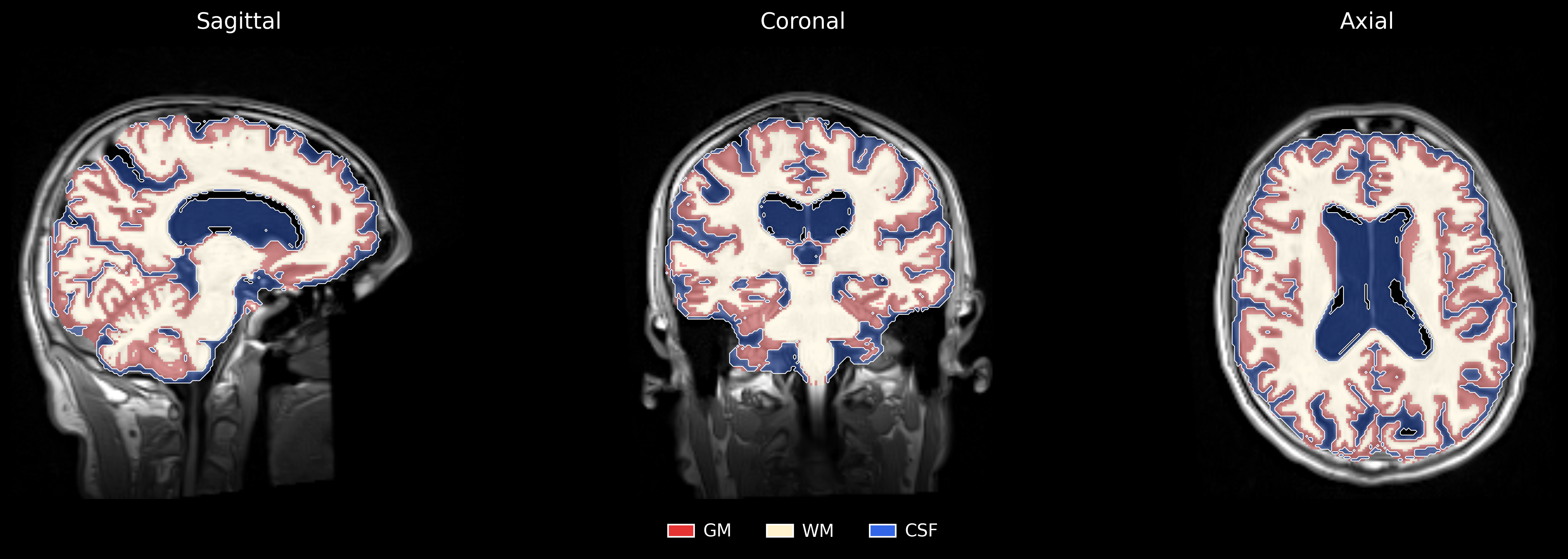}
\caption{\label{fig:qc_csf_integrity}Excluded segmentation with discontinuities in the CSF mask (outlined in white). Tissue classes: GM (red), WM (cream), CSF (blue).}
\end{figure*}

All tissue segmentations produced by sMRIPrep were visually inspected by a single rater (D.Q.) using the segmentation overlays generated during preprocessing. Segmentations were excluded when they showed gross tissue misclassification, skull-stripping failures, or discontinuities in the CSF mask (Figures \ref{fig:qc_gross_failure}, \ref{fig:qc_csf_integrity}). Such errors would introduce spurious \(H_1\) and \(H_2\) cycles in the pipeline inputs.

For the 3T cross-sectional dataset (Pipeline 1), \(1{,}161\) scans were preprocessed and \(765\) (\(65.9\%\)) passed visual inspection; \(396\) scans were excluded (\(34.1\%\)). The majority of failures involved incomplete skull stripping or boundary misclassification in temporal and inferior regions. The relatively high exclusion rate reflects the age profile of the ADNI cohort (mean age \(\approx 75\) years across groups): advanced cortical thinning and reduced tissue contrast in older brains increase the rate of automated segmentation failures, particularly in temporal and inferior regions where tissue boundaries are thinnest. Exclusion rates were comparable across diagnostic groups (CN: \(35.2\%\); MCI: \(32.6\%\); AD: \(37.3\%\)), indicating that the quality control procedure did not introduce diagnostic bias. For the 1.5T longitudinal dataset (Pipeline 2), \(170\) scans from \(85\) subjects were preprocessed; \(3\) subjects (\(6\) scans, \(3.5\%\)) were excluded, yielding the final cohort of \(82\) subjects.

\end{document}